%% file: Eta-correspondence.tex
\documentclass[a4paper]{amsart}
\usepackage{amscd, amssymb, color, booktabs, chemarrow}
\usepackage{myown}
\usepackage[makeroom]{cancel}

\usepackage{harpoon}

\title[On Theta and Eta Correspondences]{On Theta and Eta Correspondences for Finite Symplectic/Orthogonal Dual Pairs}
\author{Shu-Yen Pan}
\address{Department of Mathematics,
National Tsing Hua University, Hsinchu 300, Taiwan}
\email{sypan@math.nthu.edu.tw}


\date{\today}

\begin{document}

\begin{abstract}
\input abstract

\end{abstract}

\maketitle
\tableofcontents

\input sect01
\input sect02

\input sect03

\input sect04

\input sect05

\input sect06

\input sect07

\bibliography{refer}
\bibliographystyle{amsalpha}

\end{document}

%% file: abstract.tex
In this paper, we propose two maximal one-to-one sub-relations $\underline\theta, \overline\theta$ of
the Howe correspondence $\Theta$ for a finite reductive dual pair consisting of a symplectic group and an orthogonal group.
Moreover, we show that both the correspondences $\underline\theta$ and $\overline\theta$
are extensions beyond the stable range of the $\eta$-correspondence defined by Gurevich-Howe.

%% file: sect01.tex
% !TEX root = finite-unipotent.tex

\section{Introduction}

\subsection{}
Let $(\bfG,\bfG')$ be a reductive dual pair consisting of a symplectic group and an orthogonal group
over a finite field $\bff_q$ of odd characteristic,
and let $G$ (resp.~$G'$) denote the group of rational points of $\bfG$ (resp.~$\bfG'$).
By restricting the Weil character with respect to a non-trivial character $\psi$ of $\bff_q$
to $G\times G'$, we obtain a decomposition
\[
\omega^\psi_{\bfG,\bfG'}
=\sum_{\rho\in\cale(G),\ \rho'\in\cale(G')}m_{\rho,\rho'}\rho\otimes\rho'
\]
where $m_{\rho,\rho'}\in\bbN\cup\{0\}$,
$\cale(G)$ denotes the set of irreducible characters of $G$.
Define
\[
\Theta_{\bfG,\bfG'}=\{\,(\rho,\rho')\in\cale(G)\times\cale(G')\mid m_{\rho,\rho'}\neq 0\,\},
\]
which gives a relation from $\cale(G)$ to $\cale(G')$.
For $\rho\in\cale(G)$,
let $\Theta_{\bfG'}\colon\cale(G)\rightarrow\cale(G')$ be given by
\[
\Theta_{\bfG'}(\rho)=\{\,\rho'\in\cale(G')\mid(\rho,\rho')\in\Theta_{\bfG,\bfG'}\,\}.
\]
The subset $\{\,\rho\in\cale(G)\mid\Theta_{\bfG'}(\rho)\neq\emptyset\,\}$
is called the \emph{domain} of $\Theta_{\bfG'}$.
We say that $\rho$ \emph{occurs} in the correspondence $\Theta$ for the dual pair $(\bfG,\bfG')$ if
$\Theta_{\bfG'}(\rho)\neq\emptyset$.

It is known that
\begin{itemize}
\item $\Theta$ is \emph{symmetric}, i.e.,
$(\rho,\rho')\in\Theta_{\bfG,\bfG'}$ if and only if $(\rho',\rho)\in\Theta_{\bfG',\bfG}$.

\item $\Theta$ is \emph{persistent}, i.e.,
for $\rho\in\cale(G)$, if $\Theta_{\bfG'_{n'}}\neq\emptyset$,
then $\Theta_{\bfG'_{n''}}\neq\emptyset$ for any $n''\geq n'$ where
$\bfG'_{n'}$ denotes the group of split rank $n'$ in the Witt series of $\bfG'$.
\end{itemize}
However,
the correspondence $\Theta$ is in general not one-to-one.
More precisely, there is a dual pair $(\bfG,\bfG')$ with $\rho\in\cale(G)$ such that
$\Theta_{\bfG'}(\rho)$ contains more than one elements.

So now a natural question is: can we find a sub-relation of $\Theta$ which is one-to-one and still
has some nice properties of $\Theta$?

\subsection{}
One possible candidate for such a one-to-one sub-relation called \emph{eta correspondence} is proposed in \cite{gurevich-howe}.
Gurevich-Howe consider the dual pair $(\bfG,\bfG')=(\rmO^\epsilon_k,\Sp_{2n'})$ where $\epsilon=\pm$ and $k\leq n'$,
i.e., the dual pair is in stable range.
They show that for $\rho\in\cale(\rmO^\epsilon_k(q))$,
there is a unique $\eta(\rho)\in\Theta_{\bfG'}(\rho)\subset\cale(\Sp_{2n'}(q))$ of maximal rank in $\Theta_{\bfG'}(\rho)$.
Moreover, the mapping
\[
\eta\colon\cale(\rmO^\epsilon_k(q))\longrightarrow\cale(\Sp_{2n'}(q))
\]
is one-to-one and persistent.
The eta correspondence is only defined for a dual pair in stable range and only from the smaller group to the bigger group.
The main purpose of this paper is to propose two maximal one-to-one correspondences which extend
$\eta$ to a dual pair not necessarily in stable range.

\subsection{}
A correspondence $\underline\theta$ on unipotent characters
for certain dual pair $(\bfG,\bfG')$ is proposed in \cite{akp} under assuming the conjecture in \cite{amr}
on the $\Theta$-correspondence of unipotent characters (see also \cite{chavez}).
To describe the correspondence $\underline\theta$ we need to introduce some notations.

A \emph{$\beta$-set} is a finite subset $A=\{a_1,a_2,\ldots,a_m\}$ of non-negative integers written in
decreasing order.
A (reduced) \emph{symbol} $\Lambda=\binom{A}{B}$ is an ordered pair of two $\beta$-sets $A,B$ such that
$0\not\in A\cap B$.
For a symbol $\Lambda=\binom{a_1,a_2,\ldots,a_{m_1}}{b_1,b_2,\ldots,b_{m_2}}$,
we define a bi-partition $\Upsilon(\Lambda)$ by
\[
\Upsilon\colon \binom{a_1,\ldots,a_{m_1}}{b_1,\ldots,b_{m_2}}\mapsto
\sqbinom{a_1-(m_1-1),a_2-(m_1-2),\ldots,a_{m_1-1}-1,a_{m_1}}{b_1-(m_2-1),b_2-(m_2-2),\ldots,b_{m_2-1}-1,b_{m_2}}
\]

Let $\bfG$ be a symplectic group or an orthogonal group.
It is known that the set of irreducible unipotent characters $\cale(G)_1$ is parametrized by the set
$\cals_\bfG$ of symbols satisfying certain conditions (\cf.~\cite{lg}).
The unipotent character associated to $\Lambda$ is denoted by $\rho_\Lambda$.
For the dual pair $(\bfG,\bfG')=(\rmO^\epsilon_{2n},\Sp_{2n'})$ or $(\Sp_{2n},\rmO^\epsilon_{2n'})$,
it is known that the unipotent characters are preserved by the $\Theta$-correspondence.
We define $\underline\theta_{\bfG'}\colon\cale(G)_1\rightarrow\cale(G')_1$ by
$\underline\theta_{\bfG'}(\rho_\Lambda)=\rho_{\Lambda'}$
where $\Lambda'$ is the symbol such that (for a suitable $\tau$)
\[
\Upsilon(\Lambda')=
\begin{cases}
\Upsilon(\Lambda)^\rmt\cup\sqbinom{\tau}{-}, & \text{if $\epsilon=+$};\\
\Upsilon(\Lambda)^\rmt\cup\sqbinom{-}{\tau}, & \text{if $\epsilon=-$}
\end{cases}
\]
(\cf.~\cite{akp}, see also \cite{chavez}) where $\sqbinom{\lambda}{\mu}^\rmt=\sqbinom{\mu}{\lambda}$.
Because the conjecture by Aubert-Michel-Rouquier is proved in \cite{pan-finite-unipotent},
we see that $\underline\theta_{\bfG'}\subset\Theta_{\bfG'}$, i.e.,
$\underline\theta$ is a sub-relation of $\Theta$ on unipotent characters.
From the definition, it is obvious that $\underline\theta$ is one-to-one.

\subsection{}
The description of $\underline\theta$ (on unipotent characters) is explicit,
however, for $\Lambda\in\cals_\bfG$, $\underline\theta(\rho_\Lambda)$ is not necessarily an element of
maximal \emph{order} (\cf.~Subsection~\ref{0221}) (or maximal rank) in $\Theta_{\bfG'}(\rho_\Lambda)$.
So we propose another correspondence $\overline\theta$.
A naive way to define $\overline\theta(\rho_\Lambda)$ to be an element of maximal order in $\Theta_{\bfG'}(\rho_\Lambda)$.
But this naive definition will not guarantee $\overline\theta$ to be one-to-one in general.
So we need a modification.

Let $\cals_{n,\delta}\subset\cals_\bfG$ denote the set of symbols of rank $n$ and defect $\delta$
(\cf.~Subsection~\ref{0208}).
We know that $\Theta$ is a relation between $\cals_{n,\delta}\subset\cals_\bfG$ and
$\cals_{n',\delta'}\subset\cals_{\bfG'}$ for certain $\delta,\delta'$.
We define a linear order ``$<$'' on each $\cals_{n,\delta}$.
Then $\overline\theta_{\bfG'}(\Lambda)$ is defined inductively to be the smallest element in the set
of elements of maximal order in
\[
\Theta_{\bfG'}(\Lambda)\smallsetminus\{\,\overline\theta_{\bfG'}(\Lambda')\mid\Lambda'<\Lambda\,\}.
\]
Since $\overline\theta$ is defined inductively to guarantee injectivity,
it is usually not easy to describe $\overline\theta_{\bfG'}(\Lambda)$ explicitly.
However, we show that both $\underline\theta$ and $\overline\theta$ coincide for the following
two opposite situations (\cf.~Corollary~\ref{0503} and Proposition~\ref{0504}):
\begin{enumerate}
\item when $\cals_{n,\delta}$ and $\cals_{n',\delta'}$ are of the same size
(roughly speaking, when $\bfG$ and $\bfG'$ are of the ``similar'' sizes);

\item when $(\bfG,\bfG')$ is in stable range.
\end{enumerate}

\subsection{}
It is known that the set of irreducible characters $\cale(G)$ partitions into Lusztig series
\[
\cale(G)=\bigsqcup_{(s)\subset (G^*)^0}\cale(G)_s
\]
indexed by the conjugacy classes $(s)$ of semisimple elements in the connected component
of the dual group $G^*$ of $G$ (\cf.~\cite{lusztig-book}).

Now for $\rho\in\cale(G)_s$ we can define groups $\bfG^{(1)},\bfG^{(2)},\bfG^{(3)}$ and a modified
Lusztig correspondence
\[
\Xi_s\colon\cale(G)_s\rightarrow\cale(G^{(1)}\times G^{(2)}\times G^{(3)})_1.
\]
Write $\Xi_s(\rho)=\rho^{(1)}\otimes\rho^{(2)}\otimes\rho^{(3)}$ where $\rho^{(j)}\in\cale(G^{(j)})_1$.
From \cite{pan-Lusztig-correspondence} we know that correspondence $\Theta$ is \emph{compatible with the Lusztig correspondence},
i.e., the following diagram
\[
\begin{CD}
\rho @> \Theta_{\bfG,\bfG'} >> \rho' \\
@V \Xi_s VV @VV \Xi_{s'} V \\
\rho^{(1)}\otimes\rho^{(2)}\otimes\rho^{(3)} @> {\rm id}\otimes{\rm id}\otimes\Theta_{\bfG^{(3)},\bfG'^{(3)}}
>> \rho'^{(1)}\otimes\rho'^{(2)}\otimes\rho'^{(3)}
\end{CD}
\]
commutes.
So we can extend $\underline\theta$-correspondence outside the unipotent characters via the Lusztig correspondence:
\begin{align*}
\begin{split}
\underline\theta_{\bfG'}(\rho)=\Xi_{s'}^{-1}(\rho'^{(1)}\otimes\rho'^{(2)}\times\underline\theta_{\bfG'^{(3)}}(\rho^{(3)})),\\
\overline\theta_{\bfG'}(\rho)=\Xi_{s'}^{-1}(\rho'^{(1)}\otimes\rho'^{(2)}\times\overline\theta_{\bfG'^{(3)}}(\rho^{(3)})).
\end{split}
\end{align*}

Then we show that $\underline\theta$, $\overline\theta$ and $\eta$ coincide for a dual pair in stable range,
and hence both $\underline\theta$ and $\overline\theta$ can be regarded as extensions of the $\eta$-correspondence
(\cf.~Theorem~\ref{0404}):

\begin{thm*}
For the dual pair $(\bfG,\bfG')=(\rmO^\epsilon_k,\Sp_{2n'})$ such that $k\leq n'$,
we have
\[
\underline\theta_{\bfG'}(\rho)
=\overline\theta_{\bfG'}(\rho)
=\eta(\rho)
\]
for any $\rho\in\cale(\rmO^\epsilon_k(q))$.
\end{thm*}

\subsection{}
We say that $\vartheta$ be a \emph{sub-relation} of $\Theta$ if for each dual pair $(\bfG,\bfG')$,
$\vartheta_{\bfG,\bfG'}$ is a subset of $\Theta_{\bfG,\bfG'}$.
A sub-relation $\vartheta$ of $\Theta$ is called a \emph{theta-relation} if it is
\begin{itemize}
\item \emph{semi-persistent} on unipotent characters (\cf.~Subsections~\ref{0407}); and

\item symmetric, i.e., $\rho'=\vartheta_{\bfG'}(\rho)$ if and only if $\rho=\vartheta_\bfG(\rho')$;

\item compatible with Lusztig correspondence.
\end{itemize}
On the set of all one-to-one theta-relations,
we give a partial order by inclusion.
Then we have our second main result (\cf.~Corollary~\ref{0701}):

\begin{thm*}
Any theta-relation which properly contains $\underline\theta$ or $\overline\theta$ is not one-to-one.
\end{thm*}
Therefore, both the correspondences $\underline\theta$ and $\overline\theta$ are maximal one-to-one theta extensions
of the $\eta$-correspondence.

Although we only consider symplectic/orthogonal dual pairs in this article,
similar result which will be treated in another article also holds for unitary dual pairs.

\subsection{}
The contents of this article are as follows.
In Section 2, we provide the basic notations and preliminaries on bi-partitions and symbols which are needed for the article.
We also recall some results on the parametrization of unipotent characters of classical groups by Lusztig.
In Section 3, we recall some basic results on the theta correspondence on unipotent characters for a dual pair of a symplectic group
and an even orthogonal group.
In Section 4, we give the definitions of correspondences $\underline\theta$ and $\overline\theta$ on unipotent characters.
In Section 5, we study the relation between the correspondences $\eta$ and $\underline\theta$ on unipotent characters for
a dual pair in stable range.
In Section 6, we discuss the relation between Lusztig correspondence and the correspondences $\eta$, $\underline\theta$ and $\overline\theta$.
Then we show that both $\underline\theta$ and $\overline\theta$ are extensions of $\eta$ for dual pairs beyond stable range.
In the final section, we show that in the set of one-to-one theta relations, every element is maximal,
in particular, both $\underline\theta$ and $\overline\theta$ are maximal one-to-one theta relations.

The author would like to thank Prof. Wen-Tang Kuo for a very inspiring discussion.

%% file: sect02.tex
% !TEX root = finite-unipotent.tex

\section{Preliminaries}

\subsection{Bi-partitions and symbols}\label{0208}
For a partition $\lambda=[\lambda_1,\lambda_2\ldots,\lambda_r]$ (with $\lambda_1\geq\lambda_2\geq\cdots\geq\lambda_r\geq 0$),
we define $\|\lambda\|=\lambda_1+\cdots+\lambda_r$.
If $\|\lambda\|=n$, $\lambda$ is called a \emph{partition of $n$}.
The set of partitions of $n$ is denoted by $\calp(n)$.
For two partitions $\lambda=[\lambda_1,\ldots,\lambda_k]$ and $\lambda'=[\lambda'_1,\ldots,\lambda'_{k'}]$ of $n$,
we say that $\lambda<\lambda'$ in \emph{lexicographic order} if
there exists an index $k_0$ such that $\lambda_{k_0}<\lambda'_{k_0}$ and
$\lambda_i=\lambda'_i$ for each $i=1,\ldots,k_0-1$.
This gives a linear order on $\calp(n)$.

For a partition $\lambda=[\lambda_i]$, define its \emph{dual partition}
$\lambda^\rmT=[\lambda_j^*]$ by $\lambda_j^*=|\{\,i\mid\lambda_i\geq j\,\}|$ for $j\in\bbN$.
For $\lambda=[\lambda_1,\ldots,\lambda_r]\in\calp(n)$ and $\mu=[\mu_1,\ldots,\mu_s]$,
define their union $\lambda\cup\mu$ to be the partition in $\calp(n+m)$ with parts $\lambda_1,\ldots,\lambda_r,\mu_1,\ldots,\mu_s$.

A \emph{bi-partition of $n$} is an ordered pair of two partitions $\sqbinom{\mu}{\nu}$ such that
\[
\left\|\sqbinom{\mu}{\nu}\right\|:=\|\mu\|+\|\nu\|=n.
\]
The first (resp.~second) row of a bi-partition $\Sigma$ is denoted by $\Sigma^*$ (resp.~$\Sigma_*$), i.e.,
$\Sigma=\sqbinom{\Sigma^*}{\Sigma_*}$.
The set of bi-partitions of $n$ is denoted by $\calp_2(n)$.
For a bi-partition $\sqbinom{\mu}{\nu}$,
we define its \emph{transpose} $\sqbinom{\mu}{\nu}^\rmt=\sqbinom{\nu}{\mu}$.
We define the union of two bi-partitions by
\[
\sqbinom{\lambda}{\mu}\cup\sqbinom{\nu}{\xi}
:=\sqbinom{\lambda\cup\nu}{\mu\cup\xi}.
\]

A \emph{$\beta$-set} is a finite subset $A=\{a_1,\ldots,a_m\}$ (possibly empty) of non-negative integers
written in (strictly) decreasing order, i.e., $a_1>a_2>\cdots>a_m$.
A \emph{symbol} is an ordered pair
\[
\Lambda=\binom{A}{B}=\binom{a_1,a_2,\ldots,a_{m_1}}{b_1,b_2,\ldots,b_{m_2}}
\]
of two $\beta$-sets.
Define an equivalence relation on the set of symbols generated by the rule
\begin{equation}\label{0230}
\binom{a_1,a_2,\ldots,a_{m_1}}{b_1,b_2,\ldots,b_{m_2}}\sim
\binom{a_1+1,a_2+1,\ldots,a_{m_1}+1,0}{b_1+1,b_2+1,\ldots,b_{m_1}+1,0}.
\end{equation}
A symbol $\Lambda=\binom{A}{B}$ is \emph{reduced} if $0\not\in A\cap B$.
The first (resp.~second) row of a symbol $\Lambda$ is denoted by $\Lambda^*$ (resp.~$\Lambda_*$), i.e.,
$\Lambda=\binom{\Lambda^*}{\Lambda_*}$.
For a symbol $\Lambda=\binom{A}{B}$,
we define its \emph{transpose} $\Lambda^\rmt=\binom{B}{A}$, and its \emph{rank} and \emph{defect} by
\begin{align}\label{0202}
\begin{split}
{\rm rk}(\Lambda)
&= \sum_{a\in A}a+\sum_{b\in B}b-\left\lfloor\biggl(\frac{|A|+|B|-1}{2}\biggr)^2\right\rfloor \\
{\rm def}(\Lambda)
&= |A|-|B|
\end{split}
\end{align}
where $|X|$ denotes the number of elements in a finite set $X$.
It is easy to see that two equivalent symbols have the same rank and the same defect,
and
\begin{align}\label{0201}
\begin{split}
{\rm rk}(\Lambda)&\geq \left\lfloor\biggl(\frac{{\rm def}(\Lambda)}{2}\biggr)^2\right\rfloor, \\
\sum_{a\in A}a+\sum_{b\in B}b &={\rm rk}(\Lambda)+
\begin{cases}
m^2, & \text{if $|A|+|B|=2m+1$};\\
m(m-1), & \text{if $|A|+|B|=2m$}.
\end{cases}
\end{split}
\end{align}

\begin{lem}\label{0210}
Let
\[
\Lambda=\binom{a_1,a_2,\ldots,a_{m_1}}{b_1,b_2,\ldots,b_{m_2}}\quad\text{and}\quad
\Lambda'=\binom{a'_1,a'_2,\ldots,a'_{m'_1}}{b'_1,b'_2,\ldots,b'_{m'_2}}
\]
be two reduced symbols.
If $m_1=m_1'+1$, $m_2=m'_2+1$, $a_i>a'_i$ for $i=1,2,\ldots,m_1-1$,
and $b_j>b'_j$ for $j=1,2,\ldots,m_2-1$,
then ${\rm rk}(\Lambda)>{\rm rk}(\Lambda')$.
\end{lem}
\begin{proof}
Let
\[
\Lambda''=\binom{a'_1+1,a'_2+1,\ldots,a'_{m'_1}+1,0}{b'_1+1,b'_2+1,\ldots,b'_{m'_2}+1,0}.
\]
It is known that ${\rm rk}(\Lambda'')={\rm rk}(\Lambda')$.
Now $a_i\geq a'_i+1$ for $i=1,2,\ldots,m_1-1$ and $b_j\geq b'_j+1$ for $j=1,2,\ldots,m_2-1$.
Because now $\Lambda$ is reduced, we have $a_{m_1}>0$ or $b_{m_2}>0$.
Therefore ${\rm rk}(\Lambda)>{\rm rk}(\Lambda'')$ by the definition in (\ref{0202}).
\end{proof}

Let $\cals$ denote the set of reduced symbols,
and let $\cals_{n,\delta}$ denote the set of reduced symbols of rank $n$ and defect $\delta$.
The mapping
\begin{equation}\label{0209}
\Upsilon\colon \binom{a_1,\ldots,a_{m_1}}{b_1,\ldots,b_{m_2}}\mapsto
\sqbinom{a_1-(m_1-1),a_2-(m_1-2),\ldots,a_{m_1-1}-1,a_{m_1}}{b_1-(m_2-1),b_2-(m_2-2),\ldots,b_{m_2-1}-1,b_{m_2}}
\end{equation}
gives a bijection
\begin{equation}\label{0205}
\cals_{n,\delta}\rightarrow\begin{cases}
\calp_2(n-(\frac{\delta}{2})^2), & \text{if $\delta$ is even};\\
\calp_2(n-(\frac{\delta+1}{2})(\frac{\delta-1}{2})), & \text{if $\delta$ is odd}.
\end{cases}
\end{equation}
Modified from \cite{lg}, we define
\begin{align}\label{0220}
\begin{split}
\cals_{\Sp_{2n}} &=\{\,\Lambda\in\cals\mid{\rm rk}(\Lambda)=n,\ {\rm def}(\Lambda)\equiv 1\pmod 4\,\}; \\
\cals_{\rmO_{2n}^+} &=\{\,\Lambda\in\cals\mid{\rm rk}(\Lambda)=n,\ {\rm def}(\Lambda)\equiv 0\pmod 4\,\}; \\
\cals_{\rmO_{2n}^-} &=\{\,\Lambda\in\cals\mid{\rm rk}(\Lambda)=n,\ {\rm def}(\Lambda)\equiv 2\pmod 4\,\}.
\end{split}
\end{align}

\subsection{Unipotent characters}
Let $\bff_q$ be a finite field of $q$ elements where $q$ is a power of an odd prime $p$.
Let $\bfG$ be a classical group defined over $\bff_q$,
and let $G$ denote the group of rational points.
Let $\cale(G)=\Irr(G)$ denote the set of irreducible characters of $G$,
and let $\calv(G)$ denote the space of (complex valued) class functions on $G$.
Then $\calv(G)$ is an inner product space with $\cale(G)$ as an orthonormal basis.

The \emph{Lusztig series} $\cale(G)_s$ associated to the conjugacy class of a semisimple element $s$ 
in the connected component $(G^*)^0$ of the dual group $G^*$ of $G$ is given by
\[
\cale(G)_s=\{\,\rho\in\cale(G)\mid\langle\rho,R_{\bfT^*,s}\rangle\neq0\text{ for some $\bfT^*$ containing $s$}\,\}.
\]
Here $\bfT^*$ is a rational maximal torus in $\bfG^*$ and $s\in T^*$ (the rational points of $\bfT^*$)
and $R_{\bfT^*,s}$ is the Deligne-Lusztig virtual character associated to the pair $(\bfT^*,s)$
(\cf.~\cite{dl}).
Then we have a partition of $\cale(G)$ indexed by the conjugacy classes $(s)$:
\[
\cale(G)=\bigsqcup_{(s)\subset(G^*)^0}\cale(G)_s.
\]
Let $\calv(G)_s$ denote the subspace of $\calv(G)$ spanned by $\cale(G)_s$.
We need the following result (\cf.~\cite{DM} theorem 13.23, remark 13.24):
\begin{prop}[Lusztig]
Let $\bfG$ be a classical group.
There exists a bijection
\[
\grL_s\colon\cale(G)_s\rightarrow\cale(C_{G^*}(s))_1
\]
such that
\[
\langle\rho,R^\bfG_{\bfT^*,s}\rangle_G
=\langle\grL_s(\rho),R^{C_{\bfG^*}(s)}_{\bfT^*,1}\rangle_{C_{G^*}(s)}
\]
for any $\rho\in\cale(G)$.
Here $\langle,\rangle_G$ denotes the inner product of $\calv(G)$.
\end{prop}
Note that the mapping $\grL_s$ is usually not uniquely determined.
An element in $\cale(G)_1$ is called a \emph{unipotent character} of $G$.
The following fundamental result is from \cite{lg}:

\begin{prop}[Lusztig]
Let $\bfG$ be a symplectic group or an even orthogonal group.
There exists a bijection $\cals_\bfG\rightarrow\cale(G)_1$
where $\cals_\bfG$ is given in (\ref{0220}).
\end{prop}

Then the irreducible unipotent character of $G$ corresponding to the symbol $\Lambda\in\cals_\bfG$ is denoted by $\rho_\Lambda$.
If $\bfG$ is an orthogonal group,
it is known that $\rho_{\Lambda^\rmt}=\rho_\Lambda\cdot\sgn$.
Moreover, the trivial character $\bf 1_{\bfG}$ of $G$ is parametrized by
\[
{\bf 1}_\bfG=
\begin{cases}
\rho_{\binom{n}{-}}, & \text{if $\bfG=\Sp_{2n}$};\\
\rho_{\binom{n}{0}}, & \text{if $\bfG=\rmO^+_{2n}$};\\
\rho_{\binom{-}{n,0}}, & \text{if $\bfG=\rmO^-_{2n}$}.
\end{cases}
\]
Here $\binom{n}{-}$ means that the second row of the symbol is empty.

\subsection{The order of a symbol}\label{0221}
If $f(q)$ is a polynomial in $q$,
let $\deg(f)$ denote the degree of $f$.
For a $\beta$-set $A=\{a_1,\ldots,a_{m_1}\}$,
we define
\begin{align*}
\Delta(A,q) &=\prod_{a,a'\in A,\ a>a'}(q^a-q^{a'}), \\
\Theta(A,q^2) &=\prod_{a\in A}\left(\prod_{h=1}^a(q^{2h}-1)\right).
\end{align*}
Because $a_1>a_2>\cdots>a_{m_1}$,
we have
\[
\deg(\Delta(A,q))=\sum_{i=1}^{m_1}(m_1-i)a_i\quad\text{and}\quad
\deg(\Theta(A,q^2))=\sum_{i=1}^{m_1}a_i(a_i+1).
\]
For a symbol $\Lambda=\binom{A}{B}=\binom{a_1,a_2,\ldots,a_{m_1}}{b_1,b_2,\ldots,b_{m_2}}$,
we define
\[
\Pi(\Lambda,q)=\prod_{(a,b)\in A\times B}(q^a+q^b),
\]
and then we have
\[
\deg(\Pi(\Lambda,q))=\sum_{i=1}^{m_1}\sum_{j=1}^{m_2}\max(a_i,b_j).
\]
The following result is from \cite{lg}, theorem 8.2:

\begin{prop}[Lusztig]\label{0211}
Let $\Lambda=\binom{A}{B}=\binom{a_1,a_2,\ldots,a_{m_1}}{b_1,b_2,\ldots,b_{m_2}}\in\cals_\bfG$.
Then the dimension of $\rho_\Lambda$ is
\[
\rho_\Lambda(1)
=\frac{1}{2^c}|G|_{p'}\cdot \frac{\Delta(A,q)\Delta(B,q)}{\Theta(A,q^2)\Theta(B,q^2)}
\cdot\frac{\Pi(\Lambda,q)}{q^{\binom{m_1+m_2-2}{2}+\binom{m_1+m_2-4}{2}+\cdots}}
\]
where $c=\lfloor\frac{m_1+m_2-1}{2}\rfloor$ if $A\neq B$; $c=m_1=m_2$ if $A=B$;
$|G|_{p'}$ means the maximal factor of $|G|$ prime to $p$.
\end{prop}

The \emph{order} of $\Lambda\in\cals_\bfG$ is defined and denoted by $\ord(\Lambda)=\deg(\rho_\Lambda(1))$.
A symbol
\[
Z=\begin{cases}
\binom{z_1,z_3,\ldots,z_{m-2},z_m}{z_2,z_4,\ldots,z_{m-1}}, & \text{if $Z$ is of defect $1$};\\
\binom{z_1,z_3,\ldots,z_{m-1}}{z_2,z_4,\ldots,z_m}, & \text{if $Z$ is of defect $0$}
\end{cases}
\]
is called \emph{special} if $z_1\geq z_2\geq z_3\geq\cdots\geq z_{m-1}\geq z_m$.

\begin{lem}\label{0214}
Let $Z$ be a special symbol with entries $z_1\geq z_2\geq \cdots\geq z_m$.
Suppose that $k<l$.
Then $(l-k)+2(z_l-z_k-1)<0$.
\end{lem}
\begin{proof}
Note that two entries in the same row of a symbol must be different.
Therefore
\[
z_k-z_l\geq\begin{cases}
\frac{l-k}{2}, & \text{if $l-k$ is even;}\\
\frac{l-k-1}{2}, & \text{if $l-k$ is odd}.
\end{cases}
\]
for any $l>k$.
Then the lemma is proved immediately.
\end{proof}

\begin{lem}\label{0213}
Let $Z$ be a special symbol of rank $n$ with entries
$z_1,z_2,\ldots,z_m$ where $z_1\geq z_2\geq \cdots\geq z_m$.
Then $\ord(Z)$ is equal to
\[
\sum_{i=1}^m (m-i)z_i-\sum_{i=1}^m z_i(z_i+1)+
\begin{cases}
n(n+1)-\frac{1}{24}(m-1)(m-3)(2m-1), & \text{if $m$ odd};\\
n^2-\frac{1}{24}m(m-2)(2m-5), & \text{if $m$ even}.
\end{cases}
\]
\end{lem}
\begin{proof}
Write $Z=\binom{A}{B}$.
Because now $z_1\geq z_2\geq\cdots\geq z_m$,
from definition it is not difficult to see that
\begin{align*}
\deg(\Delta(A,q)\Delta(B,q)\Pi(Z,q))
&= \sum_{(z_i,z_j)\in Z\times Z,\ i\neq j}\max(z_i,z_j)
= \sum_{i=1}^m (m-i)z_i \\
\deg(\Theta(A,q^2)\Theta(B,q^2))
&= \sum_{i=1}^m z_i(z_i+1)
\end{align*}
Moreover,
\begin{align*}
\binom{m-2}{2}+\binom{m-4}{2}+\cdots
&= \begin{cases}
\frac{1}{24}(m-1)(m-3)(2m-1), & \text{if $m$ odd};\\
\frac{1}{24}m(m-2)(2m-5), & \text{if $m$ even},
\end{cases} \\
\deg(|G|_{p'})
&= \begin{cases}
n(n+1), & \text{if $\bfG=\Sp_{2n}$};\\
n^2, & \text{if $\bfG=\rmO^\pm_{2n}$}.
\end{cases}
\end{align*}
Note that $Z\in\cals_\bfG$ where $\bfG=\Sp_{2n}$ if $m$ is odd;
$\bfG=\rmO^\pm_{2n}$ if $m$ is even.
Then the lemma follows from Proposition~\ref{0211} immediately.
\end{proof}

For a special symbol $Z\in\cals_\bfG$,
let $\cals_Z$ denotes the set of symbols $\Lambda\in\cals_\bfG$ whose entries are exactly the same as
the entries of $Z$.

\begin{exam}
Let $Z=\binom{2,0}{1}\in\cals_{\Sp_4}$ is a special symbol of rank $2$ and defect $1$.
Then
\[
\textstyle
\cals_Z=\left\{\binom{2,0}{1},\binom{2,1}{0},\binom{1,0}{2},\binom{-}{2,1,0}\right\}\subset\cals_{\Sp_4}.
\]
\end{exam}

\begin{cor}\label{0216}
If $\Lambda\in\cals_Z$ for some special symbol $Z$,
then $\ord(\Lambda)=\ord(Z)$.
\end{cor}
\begin{proof}
From the proof of Lemma~\ref{0213},
we know that $\ord(\Lambda)$ depends only on the entries of $\Lambda$
but not on the positions of the entries.
\end{proof}

\begin{rem}
\begin{enumerate}
\item Suppose that $\bfG=\Sp_{2n}$.
If $\Lambda=\binom{n,n-1,\ldots,1,0}{n,n-1,\ldots,1}$,
then $m=2n+1$ and $z_{2i-1}=z_{2i}=n+1-i$ for $i=1,\ldots,n$.
Then
\begin{align*}
\sum_{i=1}^m (m-i)z_i &= \frac{1}{6}n(n+1)(8n+1), \\
\sum_{i=1}^m z_i(z_i+1) &= \frac{2}{3}n(n+1)(n+2), \\
\frac{1}{24}(m-1)(m-3)(2m-1) &= \frac{1}{6}n(n-1)(4n+1).
\end{align*}
So by Lemma~\ref{0213},
\begin{align*}
\ord(\Lambda)
&= \frac{1}{6}n(n+1)(8n+1)-\frac{2}{3}n(n+1)(n+2)+n(n+1)-\frac{1}{6}n(n-1)(4n+1) \\
&= n^2,
\end{align*}
i.e., $\rho_{\binom{n,n-1,\ldots,1,0}{n,n-1,\ldots,1}}$ is the Steinberg character of
$\Sp_{2n}(q)$.

\item Suppose that $\bfG=\rmO^+_{2n}$.
If $\Lambda=\binom{n,n-1,\ldots,1}{n-1,n-2,\ldots,0}$,
then $m=2n$ and $z_1=n$, $z_{2i}=z_{2i+1}=n-i$ for $i=1,\ldots,n-1$.
Then
\begin{align*}
\sum_{i=1}^m (m-i)z_i &= \frac{1}{6}n(8n^2-3n+1), \\
\sum_{i=1}^m z_i(z_i+1) &= \frac{1}{3}n(n+1)(2n+1), \\
\frac{1}{24}m(m-2)(2m-5) &= \frac{1}{6}n(n-1)(4n-5).
\end{align*}
So by Lemma~\ref{0213},
\begin{align*}
\ord(\Lambda)
&= \frac{1}{6}n(8n^2-3n+1)-\frac{1}{3}n(n+1)(2n+1)+n^2-\frac{1}{6}n(n-1)(4n-5) \\
&= n(n-1).
\end{align*}
So $\rho_{\binom{n,n-1,\ldots,1}{n-1,n-2,\ldots,0}}$ and $\rho_{\binom{n-1,n-2,\ldots,0}{n,n-1,\ldots,1}}$
are the two Steinberg characters of $\rmO^+_{2n}(q)$ for $n\geq 1$.

\item Suppose that $\bfG=\rmO^-_{2n}$.
If $\Lambda=\binom{n,n-1,\ldots,0}{n-1,n-2,\ldots,1}$,
then $\ord(\Lambda)=n(n-1)$ by (2) and Corollary~\ref{0216}.
So $\rho_\Lambda$ and $\rho_{\Lambda^\rmt}$ are the two Steinberg characters of $\rmO^-_{2n}(q)$ for $n\geq 1$.
\end{enumerate}
\end{rem}

Suppose that $z_1\geq z_2\geq\cdots\geq z_m$.
Now we consider the entries
\begin{equation}\label{0212}
z_1,\ldots,z_{k-1},z_k+1,z_{k+1},\ldots,z_{l-1},z_l-1,z_{l+1},\ldots,z_m
\end{equation}
for some $k<l$.
If $z_{k-1}<z_k+1$, then we have $z_{k-2}>z_{k-1}=z_k$.
Then the entries in (\ref{0212}) are the same as
\[
z_1,\ldots,z_{k-2},z_{k-1}+1,z_k,z_{k+1},\ldots,z_{l-1},z_l-1,z_{l+1},\ldots,z_m
\]
and we have $z_{k-2}\geq z_{k-1}+1$.
If $z_l-1<z_{l+1}$, then we have $z_l=z_{l+1}>z_{l+2}$.
Then the entries in (\ref{0212}) are the same as
\[
z_1,\ldots,z_{k-1},z_k+1,z_{k+1},\ldots,z_{l-1},z_l,z_{l+1}-1,z_{l+2},\ldots,z_m
\]
and we have $z_{l+1}-1\geq z_{l+2}$.
Therefore, without loss of generality,
we may assume that the sequence (\ref{0212}) is still monotonically decreasing .

\begin{lem}\label{0215}
Suppose that $Z$ is a symbol with entries $z_1,\ldots,z_m$ such that $z_1\geq\cdots\geq z_m$,
and $Z'$ is a symbol with entries
\[
z_1,\ldots,z_{k-1},z_k+1,z_{k+1},\ldots,z_{l-1},z_l-1,z_{l+1},\ldots,z_m
\]
such that $k<l$, $z_{k-1}\geq z_k+1$ and $z_l-1\geq z_{l+1}$.
Then $\ord(Z')<\ord(Z)$.
\end{lem}
\begin{proof}
Note that the sum of entries in $Z$ is the same as the sum of entries in $Z'$.
So by (\ref{0201}), we see that symbols $Z$ and $Z'$ are of the same rank.
Now
\begin{align*}
(z_k+1)(m-k)-(z_k+1)(z_k+2)-z_k(m-k)+z_k(z_k+1)
&= m-k-2z_k-2, \\
(z_l-1)(m-l)-(z_l-1)z_l-z_l(m-l)+z_l(z_l+1)
&= 2z_l-m+l.
\end{align*}
From Lemma~\ref{0213}, we see that
\begin{align*}
\ord(Z')-\ord(Z)
= m-k-2z_k+2z_l-m+l
= (l-k)+2(z_l-z_k-1).
\end{align*}
Hence the lemma follows from Lemma~\ref{0214} immediately.
\end{proof}

\begin{lem}\label{0218}
Let $Z,Z'$ be two special symbols of the same rank $n$.
Suppose that $Z$ has entries $z_1,\ldots,z_m$ with $z_1\geq z_2\geq\cdots\geq z_m$,
and $Z'$ has entries $z'_1,\ldots,z'_m$ with $z'_1\geq z'_2\geq\cdots\geq z'_m$.
If
\[
z_1+z_2+\cdots+z_i\geq z'_1+z'_2+\cdots+z'_i\qquad\text{for each $i=1,\ldots,m$},
\]
then $\ord(Z)\leq\ord(Z')$,
and the equality holds if and only if $z_i=z'_i$ for all $i=1,\ldots,m$.
\end{lem}
\begin{proof}
Because $Z,Z'$ have the same rank, we know that
\[
z_1+z_2+\cdots+z_m=z'_1+z'_2+\cdots+z'_m.
\]
Suppose that $Z,Z'$ are not identical.
So there are a smallest positive integer $i_0$ such that $z_{i_0}>z'_{i_0}$ and a smallest positive integer
$j_0$ such that $z_{j_0}<z'_{j_0}$.
By the assumption in the lemma, we know that $i_0<j_0$.
Let $Z''$ be a special symbol of entries $z''_1,\ldots,z''_m$ where
$z''_{i_0}=z'_{i_0}+1$, $z''_{j_0}=z'_{j_0}-1$ and $z''_k=z'_k$ for $k\neq i_0,j_0$.
Then it is clear that $z''_1\geq z''_2\geq \cdots\geq z''_m$ and $Z$ has rank $n$.
Moreover, we still have
\[
z_1+z_2+\cdots+z_i\geq z''_1+z''_2+\cdots+z''_i\qquad\text{for each $i=1,\ldots,m$}.
\]
Now $\ord(Z'')<\ord(Z')$ by Lemma~\ref{0215}.
Hence the lemma is proved by induction.
\end{proof}

\begin{lem}\label{0217}
Let
\[
\Lambda=\binom{a_1,\ldots,a_{m_1}}{b_1,\ldots,b_{m_2}}\quad\text{and}\quad
\Lambda'=\binom{a'_1,\ldots,a'_{m_1}}{b'_1,\ldots,b'_{m_2}}
\]
be two symbols of the same rank and the same defect.
\begin{itemize}
\item[(i)] Suppose that $b_j=b'_j$ for each $j$, and
\[
a_1+a_2+\cdots+a_i\geq a'_1+a'_2+\cdots+a'_i\qquad\text{for each $i=1,\ldots,m_1$},
\]
then $\ord(\Lambda)\leq\ord(\Lambda')$ and the equality holds if and only if $a_i=a'_i$ for all $i$.

\item[(ii)] Suppose that $a_i=a'_i$ for each $i$, and
\[
b_1+b_2+\cdots+b_j\geq b'_1+b'_2+\cdots+b'_j\qquad\text{for each $j=1,\ldots,m_2$},
\]
then $\ord(\Lambda)\leq\ord(\Lambda')$ and the equality holds if and only if $b_j=b'_j$ for all $j$.
\end{itemize}
\end{lem}
\begin{proof}
Let $Z,Z'$ be the special symbols with the same entries of $\Lambda,\Lambda'$ respectively.
Let $z_1,\ldots,z_{m_1+m_2}$ (resp.~$z'_1,\ldots,z'_{m_1+m_2}$) be the entries of $Z$ (resp.~$Z'$).
Clearly, the conditions in (i) imply the condition in Lemma~\ref{0218},
and we have $\ord(Z)\leq\ord(Z')$ and the equality holds if and only if $z_i=z'_i$ for all $i=1,\ldots,m_1+m_2$.
Then by Corollary~\ref{0216}, we have $\ord(\Lambda)\leq\ord(\Lambda')$
and the equality holds if and only if $a_i=a'_i$ for all $i$.

The proof of (ii) is similar and omitted.
\end{proof}

\begin{cor}\label{0219}
Let $\Lambda,\Lambda'$ be two symbols of the same rank and the same defect such that
\[
\Upsilon(\Lambda)=\sqbinom{\lambda_1,\ldots,\lambda_{m_1}}{\mu_1,\ldots,\mu_{m_2}}\quad\text{and}\quad
\Upsilon(\Lambda')=\sqbinom{\lambda'_1,\ldots,\lambda'_{m_1}}{\mu'_1,\ldots,\mu'_{m_2}}.
\]
\begin{itemize}
\item[(i)] Suppose that $\mu_j=\mu'_j$ for each $j$, and
\[
\lambda_1+\lambda_2+\cdots+\lambda_i\geq \lambda'_1+\lambda'_2+\cdots+\lambda'_i\qquad\text{for each $i=1,\ldots,m_1$},
\]
then $\ord(\Lambda)\leq\ord(\Lambda')$ and the equality holds if and only if $\lambda_i=\lambda'_i$ for all $i$.

\item[(ii)] Suppose that $\lambda_i=\lambda'_i$ for each $i$, and
\[
\mu_1+\mu_2+\cdots+\mu_j\geq \mu'_1+\mu'_2+\cdots+\mu'_j\qquad\text{for each $j=1,\ldots,m_2$},
\]
then $\ord(\Lambda)\leq\ord(\Lambda')$ and the equality holds if and only if $\mu_j=\mu'_j$ for all $j$.
\end{itemize}
\end{cor}
\begin{proof}
Write
\[
\Lambda=\binom{a_1,\ldots,a_{m_1}}{b_1,\ldots,b_{m_2}}\quad\text{and}\quad
\Lambda'=\binom{a'_1,\ldots,a'_{m_1}}{b'_1,\ldots,b'_{m_2}}.
\]
By (\ref{0209}) we have
$\lambda_i=a_i-(m_1-i)$, $\lambda'_i=a'_i-(m_1-i)$, $\mu_j=b_j-(m_2-j)$, and $\mu'_j=b'_j-(m_2-j)$ for each
$i,j$.
Then the corollary follows from Lemma~\ref{0217} immediately.
\end{proof}

%% file: sect03.tex
% !TEX root = finite-unipotent.tex

\section{Finite Theta Correspondence}

\subsection{Finite theta correspondence}
Let $(\bfG,\bfG')$ be a \emph{reductive dual pair} consisting of one orthogonal group and 
one symplectic group defined over $\bff_q$.
By restricting the \emph{Weil character} with respect to a non-trivial character $\psi$ of $\bff_q$
to $G\times G'$,
we have a decomposition
\[
\omega^\psi_{\bfG,\bfG'}=\sum_{\rho\in\cale(G),\ \rho'\in\cale(G')}m_{\rho,\rho'}\rho\otimes\rho'
\]
where $m_{\rho,\rho'}\in\bbN\cup\{0\}$ and we define
\[
\Theta_{\bfG,\bfG'} =\{\,(\rho,\rho')\in\cale(G)\times\cale(G')\mid m_{\rho,\rho'}\neq 0\,\}.
\]
So now $\Theta_{\bfG,\bfG'}$ is regarded as a relation between $\cale(G)$ and $\cale(G')$.
We say that $\rho\otimes\rho'$ occurs in the \emph{$\Theta$-correspondence} if $(\rho,\rho')\in\Theta_{\bfG,\bfG'}$,
i.e., if $m_{\rho,\rho'}\neq 0$ .
For $\rho\in\cale(G)$,
we define
\[
\Theta_{\bfG'}(\rho) =\{\,\rho'\in\cale(G')\mid (\rho,\rho')\in\Theta_{\bfG,\bfG'}\,\}.
\]
From the definition, it is obvious that $\Theta$ is \emph{symmetric}, i.e.,
$\rho'\in\Theta_{\bfG'}(\rho)$ if and only if $\rho\in\Theta_{\bfG}(\rho')$.

For a classical group $\bfG$,
let $\bfG_n$ denote the group of split rank $n$ in the \emph{Witt series} of $\bfG$, i.e,
\[
\bfG_n=\begin{cases}
\Sp_{2n}, & \text{if $\bfG$ is symplectic};\\
\rmO^+_{2n},\rmO_{2n+1},\rmO^-_{2n+2} & \text{if $\bfG$ is orthogonal}.
\end{cases}
\]
The dual pair $(\bfG,\bfG')=(\bfG_n,\bfG'_{n'})$ is said to be \emph{in stable range} if $n\leq\frac{n'}{2}$.
It is well known that $\Theta$ is \emph{persistent}, i.e.,
for $\rho\in\cale(G_n)$, if $\Theta_{\bfG'_{n'}}(\rho)\neq\emptyset$ for some $n'$,
then $\Theta_{\bfG'_{n''}}(\rho)\neq\emptyset$ for any $n''\geq n'$.
Moreover, it is known that $\Theta_{\bfG'_{n'}}(\rho)\neq\emptyset$ if $n'\geq 2n$.

\subsection{Eta correspondence in stable range}\label{0308}
Now we recall the main result from \cite{gurevich-howe}.
Let $\rho'$ be an irreducible character of $\Sp_{2n'}(q)$.
A notion of the \emph{rank} of $\rho'$ is given in \cite{gurevich-howe}.
Let $\cale(\Sp_{2n'}(q))_{(k)}$ denote the set of irreducible characters of $\Sp_{2n'}(q)$ of rank $k$.

\begin{prop}[Gurevich-Howe]\label{0301}
Consider the dual pair $(\bfG,\bfG')=(\rmO_k^\epsilon,\Sp_{2n'})$ where $k\leq n'$.
Then
\begin{enumerate}
\item[(i)] For $\rho\in\cale(\rmO^\epsilon_k(q))$,
the set $\Theta_{\bfG'}(\rho)$ contains a unique element $\eta(\rho)$ of rank $k$,
and all other elements have ranks less than $k$.

\item[(ii)] The mapping $\rho\mapsto\eta(\rho)$ gives an embedding
$\eta\colon\cale(\rmO^\epsilon_k(q))\rightarrow\cale(\Sp_{2n'}(q))_{(k)}$.
\end{enumerate}
\end{prop}
The mapping $\eta\colon\cale(\rmO^\epsilon_k(q))\rightarrow\cale(\Sp_{2n'}(q))$
is called the \emph{$\eta$-correspondence} for the dual pair $(\rmO_k^\epsilon,\Sp_{2n'})$ 
in stable range.

\begin{cor}
Consider the dual pair $(\rmO_k^\epsilon,\Sp_{2n'})$ where $k\leq n'$.
Let $\rho'\in\cale(\Sp_{2n'}(q))$.
If $\rho'=\eta(\rho)$ for some $\rho\in\cale(\rmO^\epsilon_{k}(q))$,
then $\rho'$ does not occur in the $\Theta$-correspondence for any dual
pair $(\rmO_{k'}^\epsilon,\Sp_{2n'})$ with $k'<k$.
\end{cor}
\begin{proof}
The assumption $\rho'=\eta(\rho)$ means that $\rho'$ is of rank $k$ by the above proposition.
If $\rho'$ occurs in the $\Theta$-correspondence for a dual pair $(\rmO_{k'}^\epsilon,\Sp_{2n'})$ with $k'<k$,
then the rank of $\rho'$ is less than or equal to $k'$ by the above proposition
and we get a contradiction.
\end{proof}

\subsection{Theta correspondence of unipotent characters}\label{0207}
Let $\lambda=[\lambda_1,\ldots,\lambda_k]$ and $\mu=[\mu_1,\ldots,\mu_l]$ be two partitions.
We may assume that $k=l$ by adding some $0$'s if necessary.
Then we denote
\[
\lambda\preccurlyeq\mu\quad\text{ if \ }\mu_i-1\leq\lambda_i\leq\mu_i\text{ for each }i.
\]
The following lemma is from \cite{pan-finite-unipotent} lemma 2.15:

\begin{lem}\label{0309}
Let $\lambda=[\lambda_i]$ and $\mu=[\mu_i]$ be two partitions.
Then $\lambda^\rmT\preccurlyeq\mu^\rmT$ if and only if $\mu_{i+1}\leq\lambda_i\leq\mu_i$ for each $i$
where $\lambda^\rmT,\mu^\rmT$ denote the dual partitions of $\lambda,\mu$ respectively.
\end{lem}

\begin{lem}\label{0203}
Let $A=\{a_1,\ldots,a_m\}$ and $B=\{b_1,\ldots,b_{m'}\}$ be two $\beta$-sets.
\begin{enumerate}
\item[(i)] Suppose that $m'=m$.
Then $\Upsilon(A)^\rmT\preccurlyeq\Upsilon(B)^\rmT$ if and only if $b_1\geq a_1>b_2\geq a_2>\cdots >b_m\geq a_m$.

\item[(ii)] Suppose that $m'=m+1$.
Then $\Upsilon(A)^\rmT\preccurlyeq\Upsilon(B)^\rmT$ if and only if
$b_1>a_1\geq b_2>a_2\geq \cdots \geq b_m>a_m\geq b_{m+1}$.
\end{enumerate}
\end{lem}
\begin{proof}
The lemma follows from Lemma~\ref{0309} and the definition of $\Upsilon$ in (\ref{0209}) immediately.
\end{proof}

For two symbols $\Lambda$ and $\Lambda'$,
we write $\Upsilon(\Lambda)=\sqbinom{\lambda}{\mu}$ and $\Upsilon(\Lambda')=\sqbinom{\lambda'}{\mu'}$.
Now we define several relations on the set of symbols:
\begin{align}\label{0303}
\begin{split}
\calb^+ &=\{\,(\Lambda,\Lambda')\in\cals\times\cals
\mid\mu^\rmT\preccurlyeq\lambda'^\rmT,\ \mu'^\rmT\preccurlyeq\lambda^\rmT,
\ {\rm def}(\Lambda')=-{\rm def}(\Lambda)+1\,\};\\
\calb^- &=\{\,(\Lambda,\Lambda')\in\cals\times\cals
\mid\lambda^\rmT\preccurlyeq\mu'^\rmT,\ \lambda'^\rmT\preccurlyeq\mu^\rmT,
\ {\rm def}(\Lambda')=-{\rm def}(\Lambda)-1\,\},\\
\calb_{\Sp_{2n},\rmO^\epsilon_{2n'}} &=\calb^\epsilon\cap(\cals_{\Sp_{2n}}\times\cals_{\rmO^\epsilon_{2n'}})
\end{split}
\end{align}
where $\epsilon=\pm$.

For a dual pair $(\bfG,\bfG')$ of a symplectic group and an even orthogonal group,
it is known that the unipotent characters are preserved by the $\Theta$-correspondence,
i.e., if $(\rho,\rho')\in\Theta_{\bfG,\bfG'}$, then $\rho$ is unipotent if and only if $\rho'$ is unipotent.
So let $\omega^\psi_{\bfG,\bfG',1}$ denote the unipotent part of $\omega^\psi_{\bfG,\bfG'}$, i.e.,
$\omega^\psi_{\bfG,\bfG',1}$ is the intersection of $\omega^\psi_{\bfG,\bfG'}$ with
$\calv(G)_1\otimes\calv(G')_1$.
The following proposition is from \cite{pan-finite-unipotent} theorem 3.34:

\begin{prop}\label{0204}
Let $(\bfG,\bfG')$ be a reductive dual pair of a symplectic group and an even orthogonal group.
Then we have
\[
\omega^\psi_{\bfG,\bfG',1}
=\sum_{(\Lambda,\Lambda')\in\calb_{\bfG,\bfG'}}\rho_\Lambda\otimes\rho_{\Lambda'}
\]
\end{prop}

For $\Lambda\in\cals_\bfG$,
we will denote
\[
\Theta_{\bfG'}(\Lambda)=\{\,\Lambda'\in\cals_{\bfG'}\mid(\Lambda,\Lambda')\in\calb_{\bfG,\bfG'}\,\}.
\]
Therefore by the above proposition we have $\rho_{\Lambda'}\in\Theta_{\bfG'}(\rho_\Lambda)$ if and only if
$\Lambda'\in\Theta_{\bfG'}(\Lambda)$.

Suppose that $(\Lambda,\Lambda')\in\calb_{\bfG,\bfG'}$.
From the definition in (\ref{0303}),
we can describe the relation between $\Upsilon(\Lambda)$ and $\Upsilon(\Lambda')$ as follows:
\begin{enumerate}
\item Suppose that $\epsilon=+$.
\begin{enumerate}
\item Suppose that $\bfG=\rmO^+_{2n}$ and $\bfG'=\Sp_{2n'}$.
If $\Lambda$ is of defect $4d$ for some $d\in\bbZ$,
then $\Lambda'$ is of defect $-4d+1$.
Therefore, $\Upsilon(\Lambda)$ is a bi-partition of $n-4d^2$ and $\Upsilon(\Lambda')$ is a bi-partition of
$n'-2d(2d-1)$.

\item Suppose that $\bfG=\Sp_{2n}$ and $\bfG'=\rmO^+_{2n'}$.
If $\Lambda$ is of defect $4d+1$ for some $d\in\bbZ$,
then $\Lambda'$ is of defect $-4d$.
Therefore, $\Upsilon(\Lambda)$ is a bi-partition of $n-2d(2d+1)$ and $\Upsilon(\Lambda')$ is a bi-partition of
$n'-4d^2$.
\end{enumerate}
For both cases,
we have
\[
\|\Upsilon(\Lambda')\|=\|\Upsilon(\Lambda)\|+(n'-n)+2d.
\]
Therefore the Young diagram of $\Upsilon(\Lambda')_*$ can be obtained from the Young diagram of $\Upsilon(\Lambda)^*$
by removing $k$ boxes (for some $k\geq 0$) such that no two boxes are removing from the same column;
and the Young diagram of $\Upsilon(\Lambda')^*$ can be obtained from the Young diagram of $\Upsilon(\Lambda)_*$
by adding $k+(n'-n)+2d$ boxes such that no two boxes are adding to the same column.

\item Suppose that $\epsilon=-$.
\begin{enumerate}
\item Suppose that $\bfG=\rmO^-_{2n}$ and $\bfG'=\Sp_{2n'}$.
If $\Lambda$ is of defect $4d+2$ for some $d\in\bbZ$,
then $\Lambda'$ is of defect $-4d-3$.
Therefore, $\Upsilon(\Lambda)$ is a bi-partition of $n-(2d+1)^2$ and $\Upsilon(\Lambda')$ is a bi-partition of
$n'-(2d+1)(2d+2)$.

\item Suppose that $\bfG=\Sp_{2n}$ and $\bfG'=\rmO^-_{2n'}$.
If $\Lambda$ is of defect $4d+1$ for some $d\in\bbZ$,
then $\Lambda'$ is of defect $-4d-2$.
Therefore, $\Upsilon(\Lambda)$ is a bi-partition of $n-2d(2d+1)$ and $\Upsilon(\Lambda')$ is a bi-partition of
$n'-(2d+1)^2$.
\end{enumerate}
For both cases,
we have
\[
\|\Upsilon(\Lambda')\|=\|\Upsilon(\Lambda)\|+(n'-n)-1-2d.
\]
Therefore the Young diagram of $\Upsilon(\Lambda')^*$ can be obtained from the Young diagram of $\Upsilon(\Lambda)_*$
by removing $k$ boxes (for some $k\geq 0$) such that no two boxes are removing from the same column;
and the Young diagram of $\Upsilon(\Lambda')_*$ can be obtained from the Young diagram of $\Upsilon(\Lambda)^*$
by adding $k+(n'-n)-1-2d$ boxes such that no two boxes are adding to the same column.
\end{enumerate}
For $k\geq 0$, we define
\[
\Theta_{\bfG'}(\Lambda)_k =\begin{cases}
\{\,\Lambda'\in\Theta_{\bfG'}(\Lambda)\mid\|\Upsilon(\Lambda')_*\|=\|\Upsilon(\Lambda)^*\|-k\,\}, & \text{if $\epsilon=+$};\\
\{\,\Lambda'\in\Theta_{\bfG'}(\Lambda)\mid\|\Upsilon(\Lambda')^*\|=\|\Upsilon(\Lambda)_*\|-k\,\}, & \text{if $\epsilon=-$}.
\end{cases}
\]
Clearly, as explained above we have
\begin{equation}\label{0310}
\Theta_{\bfG'}(\Lambda)=\bigsqcup_{k\geq 0}\Theta_{\bfG'}(\Lambda)_k
\end{equation}
and it is known that $\Theta_{\bfG'}(\Lambda)_k=\emptyset$ if either $\epsilon=+$ and $k>\lambda_1$; or
$\epsilon=-$ and $k>\mu_1$.

%% file: sect04.tex
% !TEX root = finite-unipotent.tex

\section{Two One-to-one Correspondences}
In this section, we consider a reductive dual pair $(\bfG,\bfG')$ of one even orthogonal group and one symplectic group.

\subsection{Definition of $\theta_k$}
Let $(\bfG,\bfG')=(\rmO^\epsilon_{2n},\Sp_{2n'})$ or $(\Sp_{2n},\rmO^\epsilon_{2n'})$
for some non-negative integers $n,n'$.
Let
\[
\Lambda=\binom{A}{B}=\binom{a_1,a_2,\ldots,a_{m_1}}{b_1,b_2,\ldots,b_{m_2}}\in\cals_\bfG,\quad\text{and}\quad
\Upsilon(\Lambda)=\sqbinom{\lambda}{\mu}=\sqbinom{\lambda_1,\lambda_2,\ldots,\lambda_{m_1}}{\mu_1,\mu_2,\ldots,\mu_{m_2}}.
\]
Suppose that
\begin{equation}\label{0431}
{\rm def}(\Lambda)=\delta=\begin{cases}
4d+1, & \text{if $\bfG=\Sp_{2n}$};\\
4d, & \text{if $\bfG=\rmO^+_{2n}$};\\
4d+2, & \text{if $\bfG=\rmO^-_{2n}$}
\end{cases}
\end{equation}
for some $d\in\bbZ$,
and define
\[
\tau=\begin{cases}
n'-n+2d, & \text{if $\epsilon=+$};\\
n'-n-1-2d, & \text{if $\epsilon=-$}.
\end{cases}
\]
So the integer $\tau$ depends on $n,n',\epsilon$ and the defect of $\Lambda$.
From the discussion in Subsection~\ref{0207},
we know that if $(\Lambda,\Lambda')\in\calb_{\bfG,\bfG'}$,
then
\[
\|\Upsilon(\Lambda')\|=\|\Upsilon(\Lambda)\|+\tau.
\]
By interchanging the roles of $\bfG$ and $\bfG'$ if necessary,
we will assume that $\tau\geq 0$.

\subsubsection{}\label{0412}
For $\epsilon=+$ and $0\leq k\leq\lambda_1$,
we define
\[
\theta_k\colon\sqbinom{\lambda_1,\ldots,\lambda_{m_1}}{\mu_1,\ldots,\mu_{m_2}}
\mapsto \sqbinom{\mu_1,\ldots,\mu_{m_2}}{\lambda_2,\ldots,\lambda_{m_1}}\cup\sqbinom{\tau+k}{\lambda_1-k}.
\]
\begin{itemize}
\item Suppose that $\bfG=\Sp_{2n}$.
Then $\theta_k$ is a mapping from
$\calp_2(n-2d(2d+1))$ to $\calp_2(n'-4d^2)$.
Moreover, $\theta_k$ induces a mapping, still denoted by $\theta_k$,
from $\cals_{n,4d+1}$ to $\cals_{n',-4d}$,
i.e., we have the following commutative diagram:
\[
\begin{CD}
\cals_{n,4d+1} @> \theta_k >> \cals_{n',-4d} \\
@V \Upsilon VV @VV \Upsilon V \\
\calp_2(n-2d(2d+1)) @> \theta_k >> \calp_2(n'-4d^2). \\
\end{CD}
\]

\item Suppose that $\bfG=\rmO^+_{2n}$.
Then $\theta_k$ is a mapping from
$\calp_2(n-4d^2))$ to $\calp_2(n'-2d(2d-1))$.
Then we have a commutative diagram:
\[
\begin{CD}
\cals_{n,4d} @> \theta_k >> \cals_{n',-4d+1} \\
@V \Upsilon VV @VV \Upsilon V \\
\calp_2(n-4d^2) @> \theta_k >> \calp_2(n'-2d(2d-1)). \\
\end{CD}
\]
\end{itemize}

It is not difficult to see that $\theta_k(\Lambda)\in\Theta_{\bfG'}(\Lambda)_k$.
Let $i,j$ be the indices such that $\mu_{j-1}>\tau+k\geq\mu_j$ and $\lambda_{i-1}\geq\lambda_1-k>\lambda_i$.
Then we see that
\[
\theta_k(\Lambda)
=\binom{b_1+1,\ldots,b_{j-1}+1,\tau+k+m_2+1-j,b_j,\ldots,b_{m_2}}
{a_2+1,\ldots,a_{i-1}+1,a_1-k+1-i,a_i,\ldots,a_{m_1}}.
\]
In particular, we have
\begin{equation}\label{0425}
\theta_0(\Lambda)
=\begin{cases}
\binom{b_1+1,\ldots,b_{m_2}+1,0}{a_1,\ldots,a_{m_1}}, & \text{if $\tau=0$};\\
\binom{\tau+m_2,b_1,\ldots,b_{m_2}}{a_1,\ldots,a_{m_1}}, & \text{if $\tau\geq\mu_1$}.
\end{cases}
\end{equation}

\subsubsection{}
For $\epsilon=-$ and $0\leq k\leq\mu_1$,
we define
\[
\theta_k\colon\sqbinom{\lambda_1,\ldots,\lambda_{m_1}}{\mu_1,\ldots,\mu_{m_2}}
\mapsto\sqbinom{\mu_2,\ldots,\mu_{m_2}}{\lambda_1,\ldots,\lambda_{m_1}}\cup\sqbinom{\mu_1-k}{\tau+k}.
\]
\begin{itemize}
\item Suppose that $\bfG=\Sp_{2n}$.
Similar to the above case, we have the following commutative diagram:
\[
\begin{CD}
\cals_{n,4d+1} @> \theta_k >> \cals_{n',-4d-2} \\
@V \Upsilon VV @VV \Upsilon V \\
\calp_2(n-2d(2d+1)) @> \theta_k >> \calp_2(n'-(2d+1)^2). \\
\end{CD}
\]

\item Suppose that $\bfG=\rmO^-_{2n}$.
Then we have the following commutative diagram:
\[
\begin{CD}
\cals_{n,4d+2} @> \theta_k >> \cals_{n',-4d-3} \\
@V \Upsilon VV @VV \Upsilon V \\
\calp_2(n-(2d+1)^2) @> \theta_k >> \calp_2(n'-(2d+1)(2d+2)). \\
\end{CD}
\]
\end{itemize}

Again, we have $\theta_k(\Lambda)\in\Theta_{\bfG'}(\Lambda)_k$.
Let $i,j$ be the indices such that $\mu_{j-1}\geq\mu_1-k>\mu_j$ and $\lambda_{i-1}>\tau+k\geq\lambda_i$.
Then we see that
\[
\theta_k(\Lambda)
=\binom{b_2+1,\ldots,b_{j-1}+1,b_1-k+1-j,b_j,\ldots,b_{m_2}}
{a_1+1,\ldots,a_{i-1}+1,\tau+k+m_1+1-i,a_i,\ldots,a_{m_1}}.
\]
In particular, we have
\begin{equation}\label{0426}
\theta_0(\Lambda)
=\begin{cases}
\binom{b_1,\ldots,b_{m_2}}{a_1+1,\ldots,a_{m_1}+1,0}, & \text{if $\tau=0$};\\
\binom{b_1,\ldots,b_{m_2}}{\tau+m_1,a_1,\ldots,a_{m_1}}, & \text{if $\tau\geq\lambda_1$}.
\end{cases}
\end{equation}

\begin{rem}
Note that $\theta_0$ is modified from the ``$\underline\theta^{N',N}$'' in \cite{akp} definition 5.
\end{rem}

\begin{rem}\label{0429}
Now $\theta_k\colon\cals_{n,\delta}\rightarrow\cals_{n',\delta'}$ where
$\delta'=-\delta+1$ if $\epsilon=+$; and $\delta'=-\delta-1$ if $\epsilon=-$.
\begin{enumerate}
\item If $\tau=0$, then $\theta_0$ is surjective,
in fact $\Lambda'=\theta_0(\Upsilon^{-1}(\Upsilon(\Lambda')^\rmt))$ for any
$\Lambda'\in\cals_{n',\delta'}$.

\item More generally, it is clear from the definition that a symbol $\Lambda'\in\cals_{n',\delta'}$ such that
$\Upsilon(\Lambda')=\sqbinom{\mu_1',\ldots,\mu'_{m'_2}}{\lambda'_1,\ldots,\lambda'_{m'_1}}$ is in the image
of $\theta_k$ if and only if
\[
\begin{cases}
\tau+k=\mu_i'\text{ for some $i$}, & \text{when $\epsilon=+$};\\
\tau+k=\lambda_i'\text{ for some $i$}, & \text{when $\epsilon=-$}.
\end{cases}
\]
\end{enumerate}
\end{rem}

\subsection{Properties of $\theta_k$}
Let $(\bfG,\bfG')=(\rmO^\epsilon_{2n},\Sp_{2n'})$ or $(\Sp_{2n},\rmO^\epsilon_{2n'})$,
and let $\Lambda\in\cals_{n,\delta}\subset\cals_\bfG$.

\begin{lem}\label{0418}
Let $\Lambda'=\theta_k(\Lambda)$ and $\Lambda''\in\Theta_{\bfG'}(\Lambda)_k$ and write
\[
\Upsilon(\Lambda')=\sqbinom{\mu'_1,\ldots,\mu'_{m_2'}}{\lambda'_1,\ldots,\lambda'_{m_1'}}\quad\text{and}\quad
\Upsilon(\Lambda'')=\sqbinom{\mu''_1,\ldots,\mu''_{m_2'}}{\lambda''_1,\ldots,\lambda''_{m_1'}}.
\]
Then
\begin{itemize}
\item[(i)] $\mu''_1+\mu''_2+\cdots+\mu''_s\geq\mu'_1+\mu'_2+\cdots+\mu'_s$ for each $s=1,\ldots,m_2'$;

\item[(ii)] $\lambda''_1+\lambda''_2+\cdots+\lambda''_t\geq\lambda'_1+\lambda'_2+\cdots+\lambda'_t$ for each
$t=1,\ldots,m_1'$.
\end{itemize}
\end{lem}
\begin{proof}
Write $\Upsilon(\Lambda)=\sqbinom{\lambda_1,\ldots,\lambda_{m_1}}{\mu_1,\ldots,\mu_{m_2}}$.

First suppose that $\epsilon=+$.
Then we have $m_1'=m_1$ and $m_2'=m_2+1$.
Because now both $\Lambda'$ and $\Lambda''$ are in $\Theta_{\bfG'}(\Lambda)_k$,
we have
\begin{align}\label{0421}
\begin{split}
\lambda_1'+\cdots+\lambda'_{m_1}
&= \lambda_1''+\cdots+\lambda''_{m_1}
=\lambda_1+\cdots+\lambda_{m_1}-k \\
\mu_1'+\cdots+\mu'_{m_2+1}
&= \mu_1''+\cdots+\mu''_{m_2+1}
=\mu_1+\cdots+\mu_{m_1}+\tau+k.
\end{split}
\end{align}
Because $\Lambda''\in\Theta_{\bfG'}(\Lambda)$,
by Lemma~\ref{0309} (see also \cite{pan-finite-unipotent} lemma 2.15),
we know that
\begin{align}\label{0419}
\begin{split}
& \mu''_1\geq\mu_1\geq\mu''_2\geq\mu_2\geq\cdots\geq\mu''_{m_2}\geq\mu_{m_2}\geq\mu''_{m_2+1}, \\
& \lambda_1\geq\lambda''_1\geq\lambda_2\geq\lambda''_2\geq\cdots\geq\lambda_{m_1}\geq\lambda''_{m_1}.
\end{split}
\end{align}
Moreover, from the definition,
we have
\begin{equation}\label{0420}
\sqbinom{\mu'_1,\ldots,\mu'_{m_2'}}{\lambda'_1,\ldots,\lambda'_{m_1'}}
=\sqbinom{\mu_1,\ldots,\mu_{j-1},\tau+k,\mu_j,\ldots,\mu_{m_2}}
{\lambda_2,\ldots,\lambda_{i-1},\lambda_1-k,\lambda_i,\ldots,\lambda_{m_1}}
\end{equation}
for some indices $i,j$.

\begin{itemize}
\item If $s\leq j-1$, then by (\ref{0419}) and (\ref{0420}) we have
\[
\mu''_1+\mu''_2+\cdots+\mu''_s
\geq\mu_1+\mu_2+\cdots+\mu_s
=\mu'_1+\mu'_2+\cdots+\mu'_s.
\]

\item If $s\geq j$, then by (\ref{0419}) and (\ref{0420}) we have
\[
\mu'_{s+1}+\cdots+\mu'_{m_2+1}
=\mu_s+\cdots+\mu_{m_2}
\geq\mu''_{s+1}+\cdots+\mu''_{m_2+1}.
\]
and hence by (\ref{0421}) we have
\begin{align*}
\mu''_1+\mu''_2+\cdots+\mu''_s
&= (\mu''_1+\mu''_2+\cdots+\mu''_{m_2+1})-(\mu''_{s+1}+\cdots+\mu''_{m_2+1}) \\
&\geq (\mu'_1+\mu'_2+\cdots+\mu'_{m_2+1})-(\mu'_{s+1}+\cdots+\mu'_{m_2+1}) \\
&= \mu'_1+\mu'_2+\cdots+\mu'_s.
\end{align*}

\item If $t\leq i-2$, then by (\ref{0419}) and (\ref{0420}) we have
\[
\lambda''_1+\lambda''_2+\cdots+\lambda''_t
\geq \lambda_2+\lambda_3+\cdots+\lambda'_{t+1}
=\lambda'_1+\lambda'_2+\cdots+\lambda'_t.
\]

\item If $t\geq i-1$, then by (\ref{0419}) and (\ref{0420}) we have
\[
\lambda'_{t+1}+\cdots+\lambda'_{m_1}
=\lambda_{t+1}+\cdots+\lambda_{m_1}
\geq \lambda''_{t+1}+\cdots+\lambda''_{m_1}.
\]
and hence by (\ref{0421}) we have
\begin{align*}
\lambda''_1+\lambda''_2+\cdots+\lambda''_t
&=(\lambda''_1+\lambda''_2+\cdots+\lambda''_{m_1})-(\lambda''_{t+1}+\lambda''_2+\cdots+\lambda''_{m_1}) \\
&\geq (\lambda'_1+\lambda'_2+\cdots+\lambda'_{m_1})-(\lambda'_{t+1}+\lambda'_2+\cdots+\lambda'_{m_1}) \\
&= \lambda'_1+\lambda'_2+\cdots+\lambda'_t.
\end{align*}
\end{itemize}
Therefore the lemma for the case $\epsilon=+$ is proved.
The proof for the case $\epsilon=-$ is similar and is omitted.
\end{proof}

\begin{lem}\label{0413}
Let $\Lambda\in\cals_\bfG$.
Then $\theta_k(\Lambda)$ is the unique element of maximal order in $\Theta_{\bfG'}(\Lambda)_k$.
\end{lem}
\begin{proof}
It is known that $\theta_k(\Lambda)$ is in $\Theta_{\bfG'}(\Lambda)_k$ by definition.

First suppose that $\epsilon=+$.
Let $\Lambda'=\theta_k(\Lambda)$, and let $\Lambda''\in\Theta_{\bfG'}(\Lambda)_k$.
Define $\Lambda'''$ to be the symbol such that
$\Upsilon(\Lambda''')_*=\Upsilon(\Lambda')_*$ and $\Upsilon(\Lambda''')^*=\Upsilon(\Lambda'')^*$.
Then it is clear that $\Lambda'''\in\Theta_{\bfG'}(\Lambda)_k$ by Lemma~\ref{0203}.
Write
\begin{align*}
\Lambda' &=\binom{b'_1,\ldots,b'_{m'_2}}{a'_1,\ldots,a'_{m'_1}}, &
\Lambda'' &=\binom{b''_1,\ldots,b''_{m'_2}}{a''_1,\ldots,a''_{m'_1}}, &
\Lambda''' &=\binom{b'''_1,\ldots,b'''_{m'_2}}{a'''_1,\ldots,a'''_{m'_1}}; \\
\Upsilon(\Lambda')&=\sqbinom{\mu'_1,\ldots,\mu'_{m'_2}}{\lambda'_1,\ldots,\lambda'_{m'_1}}, &
\Upsilon(\Lambda'')&=\sqbinom{\mu''_1,\ldots,\mu''_{m'_2}}{\lambda''_1,\ldots,\lambda''_{m'_1}}, &
\Upsilon(\Lambda''')&=\sqbinom{\mu'''_1,\ldots,\mu'''_{m'_2}}{\lambda'''_1,\ldots,\lambda'''_{m'_1}}.
\end{align*}

Now $\Upsilon(\Lambda')_*=\Upsilon(\Lambda''')_*$ means that
$\lambda'_j=\lambda'''_j$ for each $j=1,\ldots,m'_1$.
Now $\Lambda'''\in\Theta_{\bfG'}(\Lambda)_k$, by Lemma~\ref{0418}, we have
\[
\mu'''_1+\mu'''_2+\cdots+\mu'''_i\geq\mu'_1+\mu'_2+\cdots+\mu'_i
\]
for each $i=1,\ldots,m_2'$,
hence we have $\ord(\Lambda')\geq\ord(\Lambda''')$ by Corollary~\ref{0219}.
Similarly, $\Upsilon(\Lambda''')^*=\Upsilon(\Lambda'')^*$ means that
$\mu'''_i=\mu''_i$ for each $i=1,\ldots,m'_2$.
Then $\Upsilon(\Lambda')_*=\Upsilon(\Lambda''')_*$ and Lemma~\ref{0418} imply
\[
\lambda''_1+\lambda''_2+\cdots+\lambda''_j
\geq\lambda'_1+\lambda'_2+\cdots+\lambda'_j
=\lambda'''_1+\lambda'''_2+\cdots+\lambda'''_j
\]
for each $j=1,\ldots,m_1'$, we have $\ord(\Lambda''')\geq\ord(\Lambda'')$ by Corollary~\ref{0219}.
Therefore $\ord(\Lambda')\geq\ord(\Lambda'')$.

By Lemma~\ref{0217} again, we know that the equality holds if and only if $\Lambda''=\Lambda'$, i.e.,
$\theta_k(\Lambda)$ is the unique element of maximal order in $\Theta_{\bfG'}(\Lambda)_k$.

The proof for the case $\epsilon=-$ is similar and is omitted.
\end{proof}

\begin{lem}\label{0414}
There exists a unique index $k_0$ such that either
$\theta_{k_0}(\Lambda)$ or $\theta_{k_0}(\Lambda),\theta_{k_0+1}(\Lambda)$ are the elements of maximal order in the set
$\{\,\theta_k(\Lambda)\mid k\geq 0\,\}$.
\end{lem}
\begin{proof}
First suppose that $\epsilon=+$.
Write
\[
\Lambda
=\binom{a_1,a_2,\ldots,a_{m_1}}{b_1,b_2,\ldots,b_{m_2}}
\quad\text{and}\quad
\Upsilon(\Lambda)
=\sqbinom{\lambda_1,\lambda_2,\ldots,\lambda_{m_1}}{\mu_1,\mu_2,\ldots,\mu_{m_2}}.
\]
Let $0\leq k\leq\lambda_1$.
Suppose that $i,j$ are the indices such that $\mu_{j-1}>\tau+k\geq\mu_j$ and $\lambda_{i-1}\geq\lambda_1-k>\lambda_i$.
Then we have
\[
\theta_k(\Lambda)
=\binom{b_1+1,\ldots,b_{j-1}+1,\tau+k+m_2+1-j,b_j,\ldots,b_{m_2}}
{a_2+1,\ldots,a_{i-1}+1,a_1-k+1-i,a_i,\ldots,a_{m_1}}.
\]
Now we compare $\theta_k(\Lambda)$ and $\theta_{k+1}(\Lambda)$.
\begin{enumerate}
\item If $\mu_{j-1}>\tau+k+1$,
then
\[
\theta_{k+1}(\Lambda)^*=\{b_1+1,\ldots,b_{j-1}+1,\tau+k+m_2+2-j,b_j,\ldots,b_{m_2}\};
\]

\item if $\mu_{j-1}=\tau+k+1$,
then $\tau+k+m_2+1-j=b_{j-1}$ and
\[
\theta_{k+1}(\Lambda)^*=\{b_1+1,\ldots,b_{j-2}+1,\tau+k+m_2+3-j,b_{j-1},\ldots,b_{m_2}\}.
\]
\end{enumerate}
So from $\theta_k(\Lambda)$ to $\theta_{k+1}(\Lambda)$,
there is a unique entry $\alpha_k$ in the first row of $\theta_k(\Lambda)$ is changed to
$\alpha_k+1$ and all other entries in the first row is unchanged, more precisely,
$\alpha_k=\tau+k+m_2+1-j$ for case (1); and $\alpha_k=b_{j-1}+1$ for case (2).
Similarly
\begin{enumerate}
\item[(3)] if $\lambda_1-k-1>\lambda_i$,
then
\[
\theta_{k+1}(\Lambda)_*=\{a_2+1,a_3+1,\ldots,a_{i-1}+1,a_1-k-i,a_i,\ldots,a_{m_1}\};
\]

\item[(4)] if $\lambda_1-k-1=\lambda_i$,
then $a_1-k-i=a_i$ and
\[
\theta_{k+1}(\Lambda)_*=\{a_2+1,a_3+1,\ldots,a_i+1,a_1-k-1-i,a_{i+1},\ldots,a_{m_1}\}.
\]
\end{enumerate}
So from $\theta_k(\Lambda)$ to $\theta_{k+1}(\Lambda)$,
there is a unique entry $\beta_k$ in the second row of $\theta_k(\Lambda)$
is changed to $\beta_k-1$ and all other entries in the second row are unchanged,
more precisely,
$\beta_k=a_1-k-i$ for case (3); and $\beta_k=a_i$ for case (4).

We know that the sequence $\langle\beta_k\rangle$ is strictly decreasing and the sequence $\langle\alpha_k\rangle$
is strictly increasing.
Note that by Lemma~\ref{0215} if $\alpha_k+1\leq\beta_k-1$, we have $\ord(\theta_k(\Lambda))<\ord(\theta_{k+1}(\Lambda))$;
if $\alpha_k\geq\beta_k$, then $\ord(\theta_k(\Lambda))>\ord(\theta_{k+1}(\Lambda))$.
Comparing the two sequence $\langle\alpha_k\rangle$ and $\langle\beta_k\rangle$,
we have the following situations:
\begin{enumerate}
\item If $\alpha_0>\beta_0$,
then $\alpha_i\geq\beta_i$ for each $i\geq 0$ and hence $\theta_0(\Lambda)$ has maximal order.
So we have $k_0=0$.

\item If $\alpha_{\lambda_1}\leq \beta_{\lambda_1}$,
then $\alpha_i\leq\beta_i$ for each $i\geq 0$ and hence $\theta_{\lambda_1}(\Lambda)$ has maximal order.
So we have $k_0=\lambda_1$.

\item If there is an index $k_1$ such that $\alpha_{k_1}\leq\beta_{k_1}$ and $\alpha_{k_1+1}>\beta_{k_1+1}$.
Now $\alpha_{k_1+1}>\beta_{k_1+1}$ implies that
$\ord(\theta_i(\Lambda))>\ord(\theta_{i+1}(\Lambda))$ for each $i\geq k_1+1$.
Moreover, we have $\alpha_{k_1-1}+1\leq\alpha_{k_1}\leq\beta_{k_1}\leq\beta_{k_1-1}-1$,
so we have $\ord(\theta_{i-1}(\Lambda))<\ord(\theta_i(\Lambda))$ for each $i\leq k_1$.
Now we have the following three possible cases:
\begin{enumerate}
\item if $\alpha_{k_1}=\beta_{k_1}$, then $\ord(\theta_{k_1}(\Lambda))>\ord(\theta_{k_1+1}(\Lambda))$ and
we let $k_0=k_1$;

\item if $\alpha_{k_1}=\beta_{k_1}-1$, then $\alpha_{k_1}<\beta_{k_1}$ and $\alpha_{k_1}+1>\beta_{k_1}-1$,
and hence $\theta_{k_1}(\Lambda)$ and $\theta_{k_1+1}(\Lambda)$ has the same entries and hence
$\ord(\theta_{k_1}(\Lambda))=\ord(\theta_{k_1+1}(\Lambda))$
and we let $k_0=k_1$;

\item if $\alpha_{k_1}<\beta_{k_1}-1$, then $\alpha_{k_1}+1\leq\beta_{k_1}-1$ and hence
$\ord(\theta_{k_1}(\Lambda))<\ord(\theta_{k_1+1}(\Lambda))$ and
we let $k_0=k_1+1$.
\end{enumerate}
\end{enumerate}
From the definition of $k_0$ above,
we see that the orders of the sequence $\langle\theta_k(\Lambda)\rangle$ are either
\[
\ord(\theta_0(\Lambda))<\cdots<\ord(\theta_{k_0-1}(\Lambda))<\ord(\theta_{k_0}(\Lambda))
>\theta_{k_0+1}(\Lambda)>\cdots>\theta_{\lambda_1}(\Lambda),
\]
or
\begin{multline*}
\ord(\theta_0(\Lambda))<\cdots<\ord(\theta_{k_0-1}(\Lambda))<\ord(\theta_{k_0}(\Lambda)) \\
=\theta_{k_0+1}(\Lambda)>\ord(\theta_{k_0+2}(\Lambda))>\cdots>\theta_{\lambda_1}(\Lambda).
\end{multline*}
That is, for case (2.b), $\theta_{k_0}(\Lambda),\theta_{k_0+1}(\Lambda)$ are the two elements of maximal order;
and for other cases, $\theta_{k_0}(\Lambda)$ is the unique element of maximal order.

Next suppose that $\epsilon=-$.
As above, there are a sequence $\langle\alpha_k\rangle$ strictly decreasing and a strictly increasing
sequence $\langle\beta_k\rangle$ such that
\begin{enumerate}
\item[(4)] if $\beta_0>\alpha_0$,
then we let $k_0=0$;

\item[(5)] if $\beta_{\mu_1}\leq \alpha_{\mu_1}$,
then we let $k_0=\mu_1$.

\item[(6)] if there is an index $k_1$ such that $\beta_{k_1}\leq\alpha_{k_1}$ and $\beta_{k_1+1}>\alpha_{k_1+1}$,
then we have the following three possible cases:
\begin{enumerate}
\item if $\beta_{k_1}=\alpha_{k_1}$, then we let $k_0=k_1$;

\item if $\beta_{k_1}=\alpha_{k_1}-1$, then we let $k_0=k_1$;

\item if $\beta_{k_1}<\alpha_{k_1}-1$, then we let $k_0=k_1+1$.
\end{enumerate}
\end{enumerate}
For case (4.b), $\theta_{k_0}(\Lambda),\theta_{k_0+1}(\Lambda)$ are the two elements of maximal order;
and for other cases, $\theta_{k_0}(\Lambda)$ is the unique element of maximal order.
\end{proof}

\begin{cor}
There exists a unique index $k_0$ such that either
$\theta_{k_0}(\Lambda)$ or $\theta_{k_0}(\Lambda),\theta_{k_0+1}(\Lambda)$ are the elements of maximal order in
$\Theta_{\bfG'}(\Lambda)$.
\end{cor}
\begin{proof}
We know that
\[
\Theta_{\bfG'}(\Lambda)=\bigsqcup_{k\geq 0}\Theta_{\bfG'}(\Lambda)_k.
\]
Then the corollary follows from Lemma~\ref{0413} and Lemma~\ref{0414} immediately.
\end{proof}

\begin{rem}\label{0424}
\begin{enumerate}
\item Suppose that $\epsilon=+$ and $\tau\geq\mu_1$.
Then we see that $\alpha_0=\tau+m_2$ and $\beta_0=a_1=\lambda_1+m_1-1$.
From the proof of Lemma~\ref{0414}, we know that
\begin{itemize}
\item if $\alpha_0\geq\beta_0$, then $\theta_0(\Lambda)$ is the unique element of maximal order in $\Theta_{\bfG'}(\Lambda)$;

\item if $\alpha_0=\beta_0-1$, then $\theta_0(\Lambda),\theta_1(\Lambda)$ are the two elements of (the same) maximal order
in $\Theta_{\bfG'}(\Lambda)$.
\end{itemize}

\item Suppose that $\epsilon=-$ and $\tau\geq\lambda_1$.
Then we see that $\beta_0=\tau+m_1$, $\alpha_0=b_1=\mu_1+m_2-1$, and
\begin{itemize}
\item if $\beta_0\geq\alpha_0$, then $\theta_0(\Lambda)$ is the unique element of maximal order in $\Theta_{\bfG'}(\Lambda)$;

\item if $\beta_0=\alpha_0-1$, then $\theta_0(\Lambda),\theta_1(\Lambda)$ are the two elements of (the same) maximal order
in $\Theta_{\bfG'}(\Lambda)$.
\end{itemize}
\end{enumerate}
\end{rem}

\begin{exam}
Consider the dual pair $(\rmO_{30}^+,\Sp_{30})$.
Let $\Lambda=\binom{9,4,2,1}{5,4,2,0}\in\cals_{\rmO^+_{30}}$.
Then ${\rm def}(\Lambda)=0$, $\Upsilon(\Lambda)=\sqbinom{6,2,1,1}{2,2,1}$ and $\lambda_1=6$.
Hence $\tau=15-15+0=0$ and the sequences $\langle\theta_k(\Upsilon(\Lambda))\rangle$ and $\langle\theta_k(\Lambda)\rangle$
for $k=0,1,\ldots,6$ are
\begin{align*}
& \textstyle\sqbinom{2,2,1}{6,2,1,1},\sqbinom{2,2,1,1}{5,2,1,1},\sqbinom{2,2,2,1}{4,2,1,1},\sqbinom{3,2,2,1}{3,2,1,1},
\sqbinom{4,2,2,1}{2,2,1,1},\sqbinom{5,2,2,1}{2,1,1,1},\sqbinom{6,2,2,1}{2,1,1}, \\
& \textstyle\binom{6,5,3,1,0}{9,4,2,1},\binom{6,5,3,2,0}{8,4,2,1},\binom{6,5,4,2,0}{7,4,2,1},\binom{7,5,4,2,0}{6,4,2,1},
\binom{8,5,4,2,0}{5,4,2,1},\binom{9,5,4,2,0}{5,3,2,1},\binom{10,5,4,2,0}{5,3,2,0}.
\end{align*}
The sequence $\langle\alpha_k\rangle$ is $1,3,6,7,8,9,10$,
and the sequence $\langle\beta_k\rangle$ is $9,8,7,6,4,1,0$.
Now $\alpha_2=6=\beta_2-1$,
which is Case (2.b) in the proof of the previous lemma,
so we have $k_0=2$, and both $\theta_2(\Lambda)$ and $\theta_3(\Lambda)$ are of (the same) maximal order in the set
$\Theta_{\Sp_{30}}(\Lambda)$.
\end{exam}

\subsection{Definition of $\underline\theta$ and $\overline\theta$ on unipotent characters}\label{0428}
Let $(\bfG,\bfG')=(\rmO^\epsilon_{2n},\Sp_{2n'})$ or $(\Sp_{2n},\rmO^\epsilon_{2n'})$,
and let $\Lambda\in\cals_{n,\delta}\subset\cals_\bfG$.
We first assume that $\tau\geq 0$.

The definitions of $\underline\theta_{\bfG'}\colon\cals_\bfG\rightarrow\cals_{\bfG'}$ and
$\underline\theta_{\bfG'}\colon\cale(G)_1\rightarrow\cale(G')_1$ are simple,
we just define
\begin{equation}\label{0427}
\underline\theta_{\bfG'}(\Lambda)=\theta_0(\Lambda)\quad\text{and}\quad
\underline\theta_{\bfG'}(\rho_\Lambda)=\rho_{\theta_0(\Lambda)}.
\end{equation}
Then we have a relation $\underline\theta_{\bfG,\bfG'}$ between $\cale(G)_1$ and $\cale(G')_1$ by
\[
\underline\theta_{\bfG,\bfG'}
=\{\,(\rho_\Lambda,\rho_{\Lambda'})\in\cale(G)_1\times\cale(G')_1\mid\Lambda'=\underline\theta_{\bfG'}(\Lambda)\,\}.
\]

To define $\overline\theta_{\bfG'}$, we need to introduce a linear order ``$<$''on the set $\cals_{n,\delta}$ as follows:
\begin{enumerate}
\item Suppose that $\epsilon=+$.
Let $\Lambda,\Lambda'\in\cals_{n,\delta}$ where $\delta=4d+1$ (if $\bfG=\Sp_{2n}$) or $\delta=4d$ (if $\bfG=\rmO^+_{2n}$).
We define that $\Lambda<\Lambda'$ if either
\begin{itemize}
\item $\|\Upsilon(\Lambda)^*\|<\|\Upsilon(\Lambda')^*\|$; or

\item $\|\Upsilon(\Lambda)^*\|=\|\Upsilon(\Lambda')^*\|$ and $\Upsilon(\Lambda)^*<\Upsilon(\Lambda')^*$ in
lexicographic order; or

\item $\Upsilon(\Lambda)^*=\Upsilon(\Lambda')^*$ and $\Upsilon(\Lambda)_*<\Upsilon(\Lambda')_*$ in
lexicographic order.
\end{itemize}

\item Suppose that $\epsilon=-$.
Let $\Lambda,\Lambda'\in\cals_{n,\delta}$ where $\delta=4d+1$ (if $\bfG=\Sp_{2n}$) or $\delta=4d+2$ (if $\bfG=\rmO^-_{2n}$).
We define that $\Lambda<\Lambda'$ if either
\begin{itemize}
\item $\|\Upsilon(\Lambda)_*\|<\|\Upsilon(\Lambda')_*\|$; or

\item $\|\Upsilon(\Lambda)_*\|=\|\Upsilon(\Lambda')_*\|$ and $\Upsilon(\Lambda)_*<\Upsilon(\Lambda')_*$ in
lexicographic order; or

\item $\Upsilon(\Lambda)_*=\Upsilon(\Lambda')_*$ and $\Upsilon(\Lambda)^*<\Upsilon(\Lambda')^*$ in
lexicographic order.
\end{itemize}
\end{enumerate}

\begin{exam}
For the set $\cals_{4,0}\subset\cals_{\rmO^+_8}$, the linear order defined above is
\begin{multline*}
\textstyle
\binom{3,2,1,0}{4,3,2,1}<\binom{2,1,0}{4,2,1}<\binom{1,0}{3,2}<\binom{1,0}{4,1}<\binom{0}{4}<\binom{3,1,0}{3,2,1}<\binom{2,0}{3,1}<\binom{1}{3}
<\binom{2,1}{2,1}<\binom{2,1}{3,0}\\
\textstyle
<\binom{3,0}{2,1}<\binom{2}{2}<\binom{3,2,1}{3,1,0}<\binom{3,1}{2,0}<\binom{3}{1}<\binom{4,3,2,1}{3,2,1,0}<\binom{4,2,1}{2,1,0}<\binom{3,2}{1,0}
<\binom{4,1}{1,0}<\binom{4}{0}.
\end{multline*}
\end{exam}

Now we define $\overline\theta_{\bfG'}(\Lambda)$ inductively as follows.
Assume that $\overline\theta_{\bfG'}(\Lambda')$ is defined for all $\Lambda'<\Lambda$
and consider the set
\[
\Theta_{\bfG'}^\flat(\Lambda)
:=\Theta_{\bfG'}(\Lambda)\smallsetminus\{\,\overline\theta_{\bfG'}(\Lambda')\mid\Lambda'<\Lambda\,\}.
\]
We will see that $\Theta^\flat_{\bfG'}(\Lambda)$ is always non-empty (\cf.~Lemma~\ref{0423}).
Then we define $\overline\theta_{\bfG'}(\Lambda)$ to be the smallest element in the set of
elements of maximal order in $\Theta^\flat_{\bfG'}(\Lambda)$.
Then we have a mapping $\overline\theta_{\bfG'}\colon\cale(G)_1\rightarrow\cale(G')_1$ by
$\overline\theta_{\bfG'}(\rho_\Lambda)=\rho_{\overline\theta_{\bfG'}(\Lambda)}$,
and a relation $\overline\theta_{\bfG,\bfG'}$ between $\cale(G)_1$ and $\cale(G')_1$ by
\[
\overline\theta_{\bfG,\bfG'}
=\{\,(\rho_\Lambda,\rho_{\Lambda'})\in\cale(G)_1\times\cale(G')_1\mid\Lambda'=\overline\theta_{\bfG'}(\Lambda)\,\}.
\]

\begin{rem}
The definition of $\overline\theta$ depends crucially on the choice of the linear order ``$<$'' on $\cals_{n,\delta}$.
We hope that Corollary~\ref{0503} will justify our choice.
\end{rem}

\begin{lem}\label{0423}
Let $\Lambda\in\cals_{n,\delta}\subset\cals_\bfG$.
Then the set $\Theta^\flat_{\bfG'}(\Lambda)$ is always non-empty.
\end{lem}
\begin{proof}
Suppose that $\epsilon=+$.
Write $\Upsilon(\Lambda)=\sqbinom{\lambda_1,\ldots,\lambda_{m_1}}{\mu_1,\ldots,\mu_{m_2}}$.
Let $\Lambda'_0\in\Theta_{\bfG'}(\Lambda)_0$ given by
\[
\Upsilon(\Lambda'_0)=\sqbinom{\mu_1+\tau,\mu_2,\ldots,\mu_{m_2}}{\lambda_1,\ldots,\lambda_{m_1}}.
\]
Let $\Lambda'\in\cals_{n,\delta}$ and $\Lambda'<\Lambda$.
Now we have the following situations:
\begin{itemize}
\item Suppose that either $\|\Upsilon(\Lambda')^*\|<\|\Upsilon(\Lambda)^*\|$; or
$\|\Upsilon(\Lambda')^*\|=\|\Upsilon(\Lambda)^*\|$ and $\Upsilon(\Lambda')^*<\Upsilon(\Lambda)^*$,
it is clear that $\Lambda'_0\not\in\Theta_{\bfG'}(\Lambda')$ by Lemma~\ref{0309}.

\item Suppose that $\Upsilon(\Lambda')^*=\Upsilon(\Lambda)^*$ and $\Upsilon(\Lambda')_*<\Upsilon(\Lambda)_*$.
Write $\Upsilon(\Lambda')=\sqbinom{\lambda_1,\ldots,\lambda_{m_1}}{\mu'_1,\ldots,\mu'_{m'_2}}$.
The assumption $\Upsilon(\Lambda')_*<\Upsilon(\Lambda)_*$ means that there is an index $l\geq 1$ such that
$\mu'_i=\mu_i$ for $i=1,\ldots,l$ and $\mu'_l<\mu_l$.
Because now $\sum_{i=1}^{m_2}\mu_i=\sum_{i=1}^{m'_2}\mu'_i$,
we see that there is an index $k>l$ such that $\mu'_k>\mu_k$.
Now $k>1$, so by Lemma~\ref{0309}, we know that $\Lambda'_0\not\in\Theta_{\bfG'}(\Lambda')$.
\end{itemize}
So we see that $\Lambda'_0\not\in\Theta_{\bfG'}(\Lambda')$ for any $\Lambda'<\Lambda$,
and hence $\Lambda_0'$ is in $\Theta^\flat_{\bfG'}(\Lambda)$, i.e., $\Theta^\flat_{\bfG'}(\Lambda)\neq\emptyset$.

The proof for $\epsilon=-$ is similar and omitted.
\end{proof}

\begin{lem}\label{0430}
Both $\underline\theta_{\bfG'}$ and $\overline\theta_{\bfG'}$ are one-to-one mappings from
$\cals_{n,\delta}$ to $\cals_{n',\delta'}$ where $\delta'$ is given in Remark~\ref{0429}.
\end{lem}
\begin{proof}
Recall from the definition that
\[
\theta_0\colon\sqbinom{\lambda_1,\ldots,\lambda_{m_1}}{\mu_1,\ldots,\mu_{m_2}}\mapsto
\begin{cases}
\sqbinom{\mu_1,\ldots,\mu_{m_2}}{\lambda_1,\ldots,\lambda_{m_1}}\cup\sqbinom{\tau}{-}, & \text{if $\epsilon=+$};\\
\sqbinom{\mu_1,\ldots,\mu_{m_2}}{\lambda_1,\ldots,\lambda_{m_1}}\cup\sqbinom{-}{\tau}, & \text{if $\epsilon=-$}.
\end{cases}
\]
It is clear that $\underline\theta_{\bfG'}$ is a one-to-one mapping.
From the definition above, it is also obvious that $\overline\theta_{\bfG'}$ is a one-to-one.
\end{proof}

\begin{exam}
Consider the dual pair $(\bfG,\bfG')=(\rmO^+_{20},\Sp_{22})$.
It is not difficult to check that we have the following table:
\[
\begin{tabular}{c|lllll}
\toprule
$\rmO_{20}^+$ & $\Sp_{22}$ \\
$\Lambda$ & $\theta_0(\Lambda)$ & $\theta_1(\Lambda)$ & $\theta_2(\Lambda)$ & $\theta_3(\Lambda)$ & $\theta_4(\Lambda)$ \\
\midrule
$\binom{4,3}{3,2}$ & $\binom{4,3,1}{4,3}$ & $\binom{4,3,2}{4,2}$ & $\binom{5,3,2}{4,1}$ & $\binom{6,3,2}{4,0}$ \\
$\binom{5,2}{3,2}$ & $\binom{4,3,1}{5,2}$ & $\binom{4,3,2}{4,2}$ & $\binom{5,3,2}{3,2}$ & $\binom{6,3,2}{3,1}$ & $\binom{7,3,2}{3,0}$ \\
\bottomrule
\end{tabular}
\]
Then $\binom{4,3,2}{3,2}$ is the unique element of maximal order in both $\Theta_{\bfG'}(\binom{4,3}{3,2})$ and
$\Theta_{\bfG'}(\binom{5,2}{3,2})$.
Because $\binom{4,3}{3,2}<\binom{5,2}{3,2}$, we have $\binom{4,3,2}{3,2}\not\in \Theta^\flat_{\bfG'}(\binom{5,2}{3,2})$.
Finally, we can see that $\overline\theta_{\bfG'}(\binom{4,3}{3,2})=\binom{4,3,2}{4,2}$ and
$\overline\theta_{\bfG'}(\binom{5,2}{3,2})=\binom{5,3,2}{3,2}$.
\end{exam}

\begin{lem}
If $\Lambda$ is the smallest element in $\cals_{n,\delta}$,
then $\overline\theta_{\bfG'}(\Lambda)=\theta_0(\Lambda)$.
\end{lem}
\begin{proof}
First suppose $\epsilon=+$.
Because $\Lambda$ is the smallest element in $\cals_{n,\delta}$,
we have $\Upsilon(\Lambda)^*=[0]$ and hence
$\Theta^\flat_{\bfG'}(\Lambda)=\Theta_{\bfG'}(\Lambda)=\Theta_{\bfG'}(\Lambda)_0$.
Then $\theta_0(\Lambda)$ is the unique element of minimal order in $\Theta_{\bfG'}(\Lambda)$ by Lemma~\ref{0413}.
Therefore by definition, we have $\overline\theta_{\bfG'}(\Lambda)=\theta_0(\Lambda)$.
\end{proof}

\begin{lem}\label{0432}
Let $\Lambda\in\cals_{n,\delta}\subset\cals_\bfG$.
Suppose that $\Lambda'_1,\Lambda'_2$ are two elements of maximal order in $\Theta^\flat_{\bfG'}(\Lambda)$
such that $\Lambda'_i\in\Theta_{\bfG'}(\Lambda)_{k_i}$ for $i=1,2$ with $k_1<k_2$.
Then $\overline\theta_{\bfG'}(\Lambda)\neq\Lambda'_2$.
\end{lem}
\begin{proof}
Because $\Lambda'_i\in\Theta_{\bfG'}(\Lambda)_{k_i}$, we have
$\|\Upsilon(\Lambda'_i)^*\|=\|\Upsilon(\Lambda)_*\|+\tau+k_i$ by definition.
The assumption $k_1<k_2$ means that $\|\Upsilon(\Lambda'_1)^*\|<\|\Upsilon(\Lambda'_2)^*\|$ and
hence $\Lambda'_1<\Lambda'_2$.
Therefore $\overline\theta_{\bfG'}(\Lambda)\neq\Lambda'_2$ by definition.
\end{proof}

\begin{exam}\label{0422}
We consider the dual pair $(\bfG,\bfG')=(\rmO^+_8,\Sp_{10})$.
From the inequality in (\ref{0201}), we know that
\[
\cals_{\rmO^+_8}=\cals_{4,4}\cup\cals_{4,0}\cup\cals_{4,-4}\quad\text{and}\quad
\cals_{\Sp_{10}}=\cals_{5,-3}\cup\cals_{5,1}.
\]
Now by (\ref{0209}) and (\ref{0205}), we know that $\Upsilon$ establishes the following bijections
$\cals_{4,4}\simeq\calp_2(0)$;
$\cals_{4,0}\simeq\calp_2(4)$;
$\cals_{4,-4}\simeq\calp_2(0)$;
$\cals_{5,-3}\simeq\calp_2(3)$;
$\cals_{5,1}\simeq\calp_2(5)$.
Moreover, by Proposition~\ref{0204},
$\Theta$-correspondence for the pair $(\rmO^+_8,\Sp_{10})$ establishes a relation between $\cals_{4,4}$ and $\cals_{5,-3}$;
and a relation between $\cals_{4,0}$ and $\cals_{5,1}$.
The unique symbol $\binom{-}{3,2,1,0}$ in $\cals_{4,-4}$ does not occur in the relation
$\calb_{\bfG,\bfG'}$.

For $\Lambda=\binom{2}{2}\in\cals_{\rmO^+_8}$, we know that
both $\binom{3,1}{2},\binom{3,2}{1}$ are elements of maximal order in $\Theta^\flat_{\Sp_{10}}(\Lambda)$.
Moreover, $\binom{3,1}{2}\in\Theta_{\Sp_{10}}(\Lambda)_0$ and $\binom{3,2}{1}\in\Theta_{\Sp_{10}}(\Lambda)_1$.
Hence $\overline\theta_{\Sp_{10}}(\binom{2}{2})\neq\binom{3,2}{1}$ by Lemma~\ref{0432}.

The whole $\calb_{\bfG,\bfG'}$ is given by the following table.
A symbol $\Lambda'$ of maximal order in $\Theta_{\bfG'}(\Lambda)$ is superscripted by $\natural$ (Notation: ``$\Lambda'^\natural$'');
$\Lambda'\in\Theta_{\bfG'}(\Lambda)$ is overlined (Notation: ``$\overline{\Lambda'}$'') if
$\Lambda'=\overline\theta_{\bfG'}(\Lambda)$;
$\Lambda'\in\Theta_{\bfG'}(\Lambda)$ is cancelled out (Notation: ``$\bcancel{\Lambda'}$'') if
$\Lambda'\not\in\Theta^\flat_{\bfG'}(\Lambda)$.
The first element in $\Theta_{\bfG'}(\Lambda)_k$ is $\theta_k(\Lambda)$.
In particular, the first element in $\Theta_{\bfG'}(\Lambda)_0$ is $\theta_0(\Lambda)=\underline\theta_{\bfG'}(\Lambda)$.
\[
\begin{tabular}{r|lllll}
\toprule
$\rmO_8^+$ & $\Sp_{10}$ \\
$\Lambda$ & $\Theta_{\bfG'}(\Lambda)_0$ & $\Theta_{\bfG'}(\Lambda)_1$ & $\Theta_{\bfG'}(\Lambda)_2$ & $\Theta_{\bfG'}(\Lambda)_3$ & $\Theta_{\bfG'}(\Lambda)_4$ \\
\midrule
$\binom{3,2,1,0}{-}$ & $\overline{\binom{3}{3,2,1,0}}^\natural$ \\
\midrule
\midrule
$\binom{3,2,1,0}{4,3,2,1}$ & $\overline{\binom{5,4,3,2,1}{3,2,1,0}}^\natural,\binom{5,3,2,1}{2,1,0}$ \\
$\binom{2,1,0}{4,2,1}$ & $\overline{\binom{5,3,2,1}{2,1,0}}^\natural,\binom{4,3,1}{1,0},\binom{5,2,1}{1,0}$ \\
$\binom{1,0}{3,2}$ & $\overline{\binom{4,3,1}{1,0}}^\natural,\binom{4,2}{0}$ \\
$\binom{1,0}{4,1}$ & $\overline{\binom{5,2,1}{1,0}}^\natural,\binom{4,2}{0},\binom{5,1}{0}$ \\
$\binom{0}{4}$ & $\overline{\binom{5,1}{0}}^\natural,\binom{5}{-}$ \\
\midrule
$\binom{3,1,0}{3,2,1}$ & $\overline{\binom{4,3,2,1}{3,1,0}}^\natural,\binom{4,2,1}{2,0}$ & $\bcancel{\binom{5,3,2,1}{2,1,0}},\bcancel{\binom{5,2,1}{1,0}}$ \\
$\binom{2,0}{3,1}$ & $\overline{\binom{4,2,1}{2,0}}^\natural,\binom{3,2}{1},\binom{4,1}{1}$ & $\bcancel{\binom{4,3,1}{1,0}},\bcancel{\binom{5,2,1}{1,0}},\binom{4,2}{0},\bcancel{\binom{5,1}{0}}$ \\
$\binom{1}{3}$ & $\overline{\binom{4,1}{1}}^\natural,\binom{5,0}{1}$ & $\binom{4,2}{0},\bcancel{\binom{5,1}{0}},\binom{5}{-}$ \\
\midrule
$\binom{2,1}{2,1}$ & $\overline{\binom{3,2,1}{2,1}}^\natural,\binom{4,2,0}{2,1}$ & $\bcancel{\binom{4,2,1}{2,0}},\bcancel{\binom{4,1}{1}}$ \\
$\binom{2,1}{3,0}$ & $\overline{\binom{4,2,0}{2,1}}^\natural,\binom{5,1,0}{2,1}$ & $\binom{3,2}{1},\bcancel{\binom{4,1}{1}},\binom{5,0}{1}$ \\
$\binom{3,0}{2,1}$ & $\overline{\binom{3,2,1}{3,0}}^\natural,\binom{3,1}{2}$ & $\bcancel{\binom{4,2,1}{2,0}},\bcancel{\binom{4,1}{1}}$ & $\bcancel{\binom{5,2,1}{1,0}},\bcancel{\binom{5,1}{0}}$ \\
$\binom{2}{2}$ & $\overline{\binom{3,1}{2}}^\natural,\binom{4,0}{2}$ & $\binom{3,2}{1}^\natural,\bcancel{\binom{4,1}{1}},\binom{5,0}{1}$ &
$\binom{4,2}{0},\bcancel{\binom{5,1}{0}}, \binom{5}{-}$ \\
\midrule
$\binom{3,2,1}{3,1,0}$ & $\overline{\binom{4,3,1,0}{3,2,1}}^\natural,\binom{5,2,1,0}{3,2,1}$ & $\bcancel{\binom{4,2,0}{2,1}},\binom{5,1,0}{2,1}$ \\
$\binom{3,1}{2,0}$ & $\overline{\binom{3,2,0}{3,1}}^\natural,\binom{4,1,0}{3,1}$ & $\bcancel{\binom{4,2,0}{2,1}},\binom{5,1,0}{2,1},\bcancel{\binom{3,1}{2}},\binom{4,0}{2}$ &
$\bcancel{\binom{4,1}{1}},\binom{5,0}{1}$ \\
$\binom{3}{1}$ & $\overline{\binom{2,1}{3}}^\natural$ & $\bcancel{\binom{3,1}{2}}^\natural,\binom{4,0}{2}$ & $\bcancel{\binom{4,1}{1}},\binom{5,0}{1}$ & $\bcancel{\binom{5,1}{0}},\binom{5}{-}$ \\
\midrule
$\binom{4,3,2,1}{3,2,1,0}$ & $\overline{\binom{5,3,2,1,0}{4,3,2,1}}^\natural$ & $\binom{5,2,1,0}{3,2,1}$ \\
$\binom{4,2,1}{2,1,0}$ & $\overline{\binom{4,2,1,0}{4,2,1}}^\natural$ & $\binom{5,2,1,0}{3,2,1},\binom{4,1,0}{3,1}$ & $\binom{5,1,0}{2,1}$ \\
$\binom{3,2}{1,0}$ & $\overline{\binom{3,1,0}{3,2}}^\natural$ & $\binom{4,1,0}{3,1}$ & $\binom{4,0}{2}$ \\
$\binom{4,1}{1,0}$ & $\overline{\binom{3,1,0}{4,1}}^\natural$ & $\binom{4,1,0}{3,1}^\natural,\binom{3,0}{3}$ & $\binom{5,1,0}{2,1}$ & $\binom{5,0}{1}$ \\
$\binom{4}{0}$ & $\binom{2,0}{4}$ & $\overline{\binom{3,0}{3}}^\natural$ & $\binom{4,0}{2}$ & $\binom{5,0}{1}$ & $\binom{5}{-}$ \\
\midrule
\midrule
$\binom{-}{3,2,1,0}$ & \\
\bottomrule
\end{tabular}
\]
From the above table, for most symbols $\Lambda\in\cals_{4,0}$,
we have $\underline\theta_{\Sp_{10}}(\Lambda)=\overline\theta_{\Sp_{10}}(\Lambda)$, the only exception is
$\underline\theta_{\Sp_{10}}(\binom{4}{0})=\binom{2,0}{4}$ and
$\overline\theta_{\Sp_{10}}(\binom{4}{0})=\binom{3,0}{3}$.
\end{exam}

Now $\underline\theta_{\bfG'}(\Lambda)$ and $\overline\theta_{\bfG'}(\Lambda)$ are defined if $\tau\geq 0$.
For case $\tau<0$, we can extend the definitions of $\underline\theta$ and $\overline\theta$ by symmetry, i.e., we define
\begin{align*}
\underline\theta_{\bfG'}(\Lambda) &= \Lambda'\quad\text{if and only if}\quad\underline\theta_\bfG(\Lambda')=\Lambda;\\
\overline\theta_{\bfG'}(\Lambda) &= \Lambda'\quad\text{if and only if}\quad\overline\theta_\bfG(\Lambda')=\Lambda.
\end{align*}
From the definition, we see that
if $\underline\theta_{\bfG'}(\Lambda)=\Lambda'$, $\Upsilon(\Lambda)=\sqbinom{\lambda_1,\ldots,\lambda_{m_1}}{\mu_1,\ldots,\mu_{m_2}}$
and $\tau<0$,
then
\[
\Upsilon(\Lambda')=\begin{cases}
\sqbinom{\mu_1,\ldots,\mu_{m_2}}{\lambda_1,\ldots,\lambda_{m_1}}\smallsetminus\sqbinom{-}{\lambda_i}, & \text{if $\epsilon=+$};\\
\sqbinom{\mu_1,\ldots,\mu_{m_2}}{\lambda_1,\ldots,\lambda_{m_1}}\smallsetminus\sqbinom{\mu_i}{-}, & \text{if $\epsilon=-$}
\end{cases}
\]
for some $i$.
So it is possible that the domain of $\underline\theta_{\bfG'}$ is properly contained in the domain of $\Theta_{\bfG'}$, i.e.,
there exists $\rho\in\cale(G)$ such that $\underline\theta_{\bfG'}(\rho)$ is not defined but
$\Theta_{\bfG'}(\rho)\neq\emptyset$.

\begin{exam}\label{0319}
Suppose that $\Lambda=\binom{4,1}{3,1}\in\cals_{\rmO_{14}^+}$.
Then $\Upsilon(\Lambda)=\sqbinom{3,1}{2,1}\in\calp_2(7)$.
By definition and Proposition~\ref{0204}, it is not difficult to see that
\begin{align*}
\Theta_{\Sp_8}(\rho_\Lambda)
&=\{\underline\theta_{\Sp_8}(\rho_\Lambda)\}=\{\rho_{\binom{3,1}{1}}\};\\
\Theta_{\Sp_{10}}(\rho_\Lambda)
&=\{\rho_{\binom{4,1}{1}},\rho_{\binom{3,2}{1}},\rho_{\binom{4,2,1}{2,0}},\rho_{\binom{4,2,0}{2,1}},\}
\end{align*}
From the above we know that $\underline\theta_{\Sp_{2n'}}(\rho_\Lambda)$ is defined only for $n'=4,6,7,8,9,\ldots$,
in particular, $\underline\theta_{\Sp_{10}}(\rho_\Lambda)$ does not exist.
\end{exam}

%% file: sect05.tex
% !TEX root = finite-unipotent.tex

\section{Theta and Eta Correspondences on Unipotent Characters}
In this section, we consider a reductive dual pair $(\bfG,\bfG')$ of one even orthogonal group and one symplectic group.

\subsection{The case for $\tau=0$}
In this subsection, we consider the situation that $\tau=0$, i.e.,
two sets $\cals_{n,\delta}$ and $\cals_{n',\delta'}$ are of the same size.

\begin{lem}\label{0505}
Suppose that $\tau=0$.
Let $\Lambda\in\cals_{n,\delta}\subset\cals_\bfG$.
Then $\Theta_{\bfG'}(\Lambda)_0=\{\theta_0(\Lambda)\}$.
\end{lem}
\begin{proof}
First suppose that $\epsilon=+$.
Let $\Lambda'\in\Theta_{\bfG'}(\Lambda)_0$.
Because now $k=\tau=0$,
we see from Subsection~\ref{0207} that $\Upsilon(\Lambda')$ is obtained from $\Upsilon(\Lambda)^\rmt$ by removing zero box from the Young diagram
of $\Upsilon(\Lambda)_*$ and adding zero box to the Young diagram of $\Upsilon(\Lambda)^*$.
This means that $\Upsilon(\Lambda')=\Upsilon(\Lambda)^\rmt$, i.e., $\Lambda'=\theta_0(\Lambda)$.

The proof for $\epsilon=-$ is similar.
\end{proof}

\begin{lem}\label{0415}
Suppose that $\tau=0$.
Let $\Lambda\in\cals_{n,\delta}\subset\cals_\bfG$.
Then $\Theta^\flat_{\bfG'}(\Lambda)=\{\theta_0(\Lambda)\}$.
\end{lem}
\begin{proof}
First suppose that $\epsilon=+$.
We now prove the lemma by induction.
\begin{itemize}
\item Suppose that $\Lambda$ is the smallest element in $\cals_{n,\delta}$.
Then $\Upsilon(\Lambda)^*=[0]$, and we know that $\Theta_{\bfG'}(\Lambda)=\Theta_{\bfG'}(\Lambda)_0$ by (\ref{0310}).
Because now we also have $\tau=0$, we have $\Theta_{\bfG'}(\Lambda)_0=\{\theta_0(\Lambda)\}$ by Lemma~\ref{0505}.
Therefore,
we have
\[
\Theta_{\bfG'}^\flat(\Lambda)
=\Theta_{\bfG'}(\Lambda)
=\Theta_{\bfG'}(\Lambda)_0
=\{\theta_0(\Lambda)\},
\]
i.e., the assertion is true for this case.

\item Now we assume that $\Theta^\flat_{\bfG'}(\Lambda')=\{\theta_0(\Lambda')\}$ for any
$\Lambda'$ such that $\Lambda'<\Lambda$.
Then by definition, we have $\overline\theta_{\bfG'}(\Lambda')=\theta_0(\Lambda')$.
Now suppose that $\Lambda''\in\Theta_{\bfG'}(\Lambda)_k$ for some $k>0$.
Let $\Lambda'''=\Upsilon^{-1}((\Upsilon(\Lambda''))^\rmt)\in\cals_{n,\delta}$.
Then $\|\Upsilon(\Lambda''')^*\|<\|\Upsilon(\Lambda)^*\|$,
in particular $\Lambda'''<\Lambda$, 
and hence $\overline\theta_{\bfG'}(\Lambda''')=\theta_0(\Lambda''')=\Lambda''$
by induction hypothesis.
This means that $\Lambda''\not\in\Theta_{\bfG'}^\flat(\Lambda)$ and hence
\[
\Theta_{\bfG'}^\flat(\Lambda)\subseteq\Theta_{\bfG'}(\Lambda)_0=\{\theta_0(\Lambda)\}
\]
by Lemma~\ref{0505}.
But we know that $\Theta_{\bfG'}^\flat(\Lambda)\neq\emptyset$ by Lemma~\ref{0423}, so we conclude
$\Theta_{\bfG'}^\flat(\Lambda)=\{\theta_0(\Lambda)\}$.
\end{itemize}
So the lemma is proved by induction for the case $\epsilon=+$.

The proof for $\epsilon=-$ is similar.
\end{proof}

\begin{cor}\label{0503}
Suppose that $\tau=0$.
Then the mapping $\overline\theta_{\bfG'}\colon\cals_{n,\delta}\rightarrow\cals_{n',\delta'}$ is given by
\[
\binom{a_1,a_2,\ldots,a_{m_1}}{b_1,b_2,\ldots,b_{m_2}}\mapsto\begin{cases}
\binom{b_1+1,b_2+1,\ldots,b_{m_2}+1,0}{a_1,a_2,\ldots,a_{m_1}}, & \text{if $\epsilon=+$};\\
\binom{b_1,b_2,\ldots,b_{m_2}}{a_1+1,a_2+1,\ldots,a_{m_1}+1,0}, & \text{if $\epsilon=-$},
\end{cases}
\]
i.e., $\underline\theta_{\bfG'}(\Lambda)=\overline\theta_{\bfG'}(\Lambda)$.
\end{cor}
\begin{proof}
Suppose that $\tau=0$ and $\Lambda\in\cals_{n,\delta}\subset\cals_\bfG$.
Then by Lemma~\ref{0415} and the definitions of $\underline\theta$ and $\overline\theta$, we
have $\overline\theta_{\bfG'}(\Lambda)=\theta_0(\Lambda)=\underline\theta_{\bfG'}(\Lambda)$.
\end{proof}

\subsection{The case in stable range}

\begin{lem}\label{0416}
Let $(\bfG,\bfG')$ be a dual pair in the stable range and $\Lambda\in\cals_\bfG$.
Then $\theta_0(\Lambda)$ is the unique element of maximal order in the set $\Theta_{\bfG'}(\Lambda)$.
\end{lem}
\begin{proof}
Write
\[
\Lambda=\binom{a_1,\ldots,a_{m_1}}{b_1,\ldots,b_{m_2}}\in\cals_\bfG\quad\text{and}\quad
\Upsilon(\Lambda)=\sqbinom{\lambda_1,\ldots,\lambda_{m_1}}{\mu_1,\ldots,\mu_{m_2}}.
\]

\begin{enumerate}
\item Suppose that $(\bfG,\bfG')=(\Sp_{2n},\rmO^+_{2n'})$ with $n'\geq 2n$, and ${\rm def}(\Lambda)=m_1-m_2=4d+1$
for some $d\in\bbZ$.
Then $\Upsilon(\Lambda)\in\calp_2(n-2d(2d+1))$ by (\ref{0205}), and
\begin{align*}
\tau-\mu_1
\geq \tau-(n-2d(2d+1))
&=n'-n+2d-(n-2d(2d+1)) \\
&=4d^2+4d\geq 0
\end{align*}
So now $\alpha_0=\tau+m_2$, $\beta_0=a_1=\lambda_1+m_1-1$,
and hence
\begin{align*}
\alpha_0-\beta_0
=\tau-\lambda_1+m_2-m_1+1
&\geq \tau-(n-2d(2d+1))+m_2-m_1+1 \\
&\geq 4d^2\geq 0.
\end{align*}

\item Suppose that $(\bfG,\bfG')=(\rmO^+_{2n},\Sp_{2n'})$ with $n'\geq 2n$,
and ${\rm def}(\Lambda)=m_1-m_2=4d$ for some $d\in\bbZ$.
Then $\Upsilon(\Lambda)\in\calp_2(n-4d^2)$ and
\[
\tau-\mu_1
\geq\tau-(n-4d^2)
=n'-n+2d-(n-4d^2)=4d^2+2d\geq 0
\]
So now $\alpha_0=\tau+m_2$, $\beta_0=a_1=\lambda_1+m_1-1$, and hence
\begin{align*}
\alpha_0-\beta_0
=\tau-\lambda_1+m_2-m_1+1
&\geq \tau-(n-4d^2)+m_2-m_1+1 \\
&\geq 4d^2-2d+1\geq 0.
\end{align*}

\item Suppose that $(\bfG,\bfG')=(\Sp_{2n},\rmO^-_{2n'})$ with $n'-1\geq 2n$,
and ${\rm def}(\Lambda)=m_1-m_2=4d+1$ for some $d\in\bbZ$.
Then $\Upsilon(\Lambda)\in\calp_2(n-2d(2d+1))$ and
\begin{align*}
\tau-\lambda_1
\geq\tau-(n-2d(2d+1))
&=n'-n-1-2d-(n-2d(2d+1)) \\
&=4d^2\geq 0
\end{align*}
So now $\beta_0=\tau+m_1$, $\alpha_0=b_1=\mu_1+m_2-1$,
and hence
\begin{align*}
\beta_0-\alpha_0
=\tau-\mu_1+m_1-m_2+1
&\geq \tau-(n-2d(2d+1))+m_1-m_2+1 \\
&\geq 4d^2+4d+2\geq 0.
\end{align*}

\item Suppose that $(\bfG,\bfG')=(\rmO^-_{2n},\Sp_{2n'})$ with $n'\geq 2n$,
and ${\rm def}(\Lambda)=m_1-m_2=4d+2$ for some $d\in\bbZ$.
Then $\Upsilon(\Lambda)\in\calp_2(n-(2d+1)^2)$ and
\begin{align*}
\tau-\lambda_1
\geq\tau-(n-(2d+1)^2)
&=n'-n-1-2d-(n-(2d+1)^2) \\
&=4d^2+2d\geq 0
\end{align*}
So now $\beta_0=\tau+m_1$, $\alpha_0=b_1=\mu_1+m_2-1$, and hence
\begin{align*}
\beta_0-\alpha_0
=\tau-\mu_1+m_1-m_2+1
&\geq \tau-(n-(2d+1)^2)+m_1-m_2+1 \\
&\geq 4d^2+6d+3\geq 0.
\end{align*}
\end{enumerate}
Now we have $\alpha_0\geq\beta_0$ when $\epsilon=+$;
and $\beta_0\geq\alpha_0$ when $\epsilon=-$,
so by Remark~\ref{0424}, we conclude that $\theta_0(\Lambda)$ is the unique element of maximal order in the set
$\Theta_{\bfG'}(\Lambda)$.
\end{proof}

\begin{lem}\label{0506}
Let $(\bfG,\bfG')$ be a dual pair in the stable range and $\Lambda\in\cals_\bfG$.
Then $\theta_0(\Lambda)\in\Theta^\flat_{\bfG'}(\Lambda)$.
\end{lem}
\begin{proof}
We prove the proposition by induction.
First suppose that $\Lambda$ is the smallest element in $\cals_{n,\delta}$.
Then the proposition clearly follows from Lemma~\ref{0416} since now $\Theta^\flat_{\bfG'}(\Lambda)=\Theta_{\bfG'}(\Lambda)$.
Now we assume that $\theta_0(\Lambda')\in\Theta^\flat_{\bfG'}(\Lambda')$ for each $\Lambda'<\Lambda$.
Lemma~\ref{0416} implies $\overline\theta_{\bfG'}(\Lambda')=\theta_0(\Lambda')$ for each $\Lambda'<\Lambda$.
Because the mapping $\theta_0$ is one-to-one by Lemma~\ref{0430}, 
we have $\theta_0(\Lambda)\in\Theta^\flat_{\bfG'}(\Lambda)$.
\end{proof}

\begin{prop}\label{0504}
Let $(\bfG,\bfG')$ be a dual pair in the stable range and $\Lambda\in\cals_\bfG$.
Then $\underline\theta_{\bfG'}(\Lambda)=\overline\theta_{\bfG'}(\Lambda)$.
\end{prop}
\begin{proof}
The proposition follows from Lemma~\ref{0416} and Lemma~\ref{0506} immediately.
\end{proof}

\subsection{The first occurrence for unipotent characters}\label{0502}
We know that both $\underline\theta_{\bfG'}(\Lambda)$ and $\overline\theta_{\bfG'}(\Lambda)$ are defined
whenever $\tau\geq 0$.
Because now $\underline\theta$ and $\overline\theta$ are symmetric,
we conclude that every irreducible unipotent character $\rho_{\Lambda'}\in\cale(G')_1$ eventually occurs in the
the correspondence $\underline\theta$ (resp.~$\overline\theta$), i.e.,
there exist a group $\bfG$ and an irreducible character $\rho_\Lambda\in\cale(G)_1$ such that
$\underline\theta_{\bfG'}(\rho_\Lambda)=\rho_{\Lambda'}$ (resp.~$\overline\theta_{\bfG'}(\rho_\Lambda)=\rho_{\Lambda'}$).

For $\rho_{\Lambda'}\in\cale(G')_1$, define $n_0(\rho_{\Lambda'})$ to be the smallest $n$
such that $\Theta_{\bfG_n}(\rho_{\Lambda'})$ is non-empty.
Similarly, define $\underline n_0(\rho_{\Lambda'})$ (resp.\ $\overline n_0(\rho_{\Lambda'})$) to be the smallest $n$
such that $\underline\theta_{\bfG_n}(\rho_{\Lambda'})$ (resp.\ $\overline\theta_{\bfG_n}(\rho_{\Lambda'})$) is defined.
Because $\underline\theta_\bfG(\Lambda')\in\Theta_\bfG(\Lambda')$ and $\overline\theta_\bfG(\Lambda')\in\Theta_\bfG(\Lambda')$,
we have $n_0(\rho_{\Lambda'})\leq\underline n_0(\rho_{\Lambda'})$ and $n_0(\rho_{\Lambda'})\leq\overline n_0(\rho_{\Lambda'})$
for any $\Lambda'\in\cals_{\bfG'}$.

\begin{lem}\label{0314}
Let $\Lambda'\in\cals_{\bfG'}$.
Then $n_0(\rho_{\Lambda'})=\underline n_0(\rho_{\Lambda'})$.
\end{lem}
\begin{proof}
Write $\Lambda'=\binom{a'_1,\ldots,a'_{m'_1}}{b'_1,\ldots,b'_{m'_2}}$ and define
\[
\Lambda=\begin{cases}
\binom{b'_1+1,\ldots,b'_{m'_2}+1,0}{-}, & \text{if $\epsilon=+$ and $m'_1=0$};\\
\binom{b'_1,\ldots,b'_{m'_2}}{a'_2,\ldots,a'_{m'_1}}, & \text{if $\epsilon=+$ and $m'_1\geq 1$};\\
\binom{-}{a'_1+1,\ldots,a'_{m'_1}+1,0}, & \text{if $\epsilon=-$ and $m'_2=0$};\\
\binom{b'_2,\ldots,b'_{m'_2}}{a'_1,\ldots,a'_{m'_1}}, & \text{if $\epsilon=-$ and $m'_2\geq 1$}.
\end{cases}
\]
Then by (\ref{0425}) and (\ref{0426}),
we have $\Lambda'=\theta_0(\Lambda)$,
and hence $\rho_{\Lambda'}=\underline\theta_{\bfG'}(\rho_\Lambda)$.
From the results in section 8 of \cite{pan-Lusztig-correspondence},
we conclude that $n_0(\rho_{\Lambda'})=\underline n_0(\rho_{\Lambda'})$.
\end{proof}

\begin{exam}
Let $\Lambda'=\binom{2,0}{4}\in\cals_{\Sp_{10}}$.
\begin{enumerate}
\item Suppose that $\epsilon=+$.
From Example~\ref{0422}, we know that $\underline\theta_{\Sp_{10}}(\binom{4}{0})=\Lambda'$ and
$\overline\theta_{\Sp_{10}}(\binom{4}{0})=\binom{3,0}{3}\neq\Lambda'$.
Therefore, we have $\underline n_0(\rho_{\Lambda'})=4$ and $\overline n_0(\rho_{\Lambda'})=5$.

\item Suppose that $\epsilon=-$.
Note that $\binom{-}{2,0}\in\cals_{\rmO^-_4}$ and $(\rmO^-_4,\Sp_{10})$ is in stable range.
Then we have $\underline\theta_{\Sp_{10}}(\binom{-}{2,0})=\overline\theta_{\Sp_{10}}(\binom{-}{2,0})=\Lambda'$,
and hence $\underline n_0(\rho_{\Lambda'})=\overline n_0(\rho_{\Lambda'})=2$.
\end{enumerate}
\end{exam}

\subsection{The correspondences $\eta$, $\underline\theta$ and $\overline\theta$ on unipotent characters}
Now we are going to show that the three correspondences $\eta$, $\underline\theta$ and $\overline\theta$
coincide on their common domain,
i.e., the situation that the dual pair $(\bfG,\bfG')$ is in stable range and the irreducible characters are unipotent.

\begin{lem}\label{0302}
Let $(\bfG,\bfG')=(\bfG_n,\bfG'_{n'})$ be a dual pair in stable range and $\Lambda\in\cals_{\bfG_n}$.
Suppose that $\Lambda'\in\Theta_{\bfG'}(\Lambda)\smallsetminus\{\underline\theta_{\bfG'}(\Lambda)\}$.
Then there exist $n''<n$ and $\Lambda_0\in\cals_{\bfG_{n''}}$ such that $\underline\theta_{\bfG'}(\Lambda_0)=\Lambda'$,
in particular, $n_0(\rho_{\Lambda'})<n$.
\end{lem}
\begin{proof}
Write
\[
\Lambda=\binom{a_1,\ldots,a_{m_1}}{b_1,\ldots,b_{m_2}},\quad
\Upsilon(\Lambda)=\sqbinom{\lambda_1,\ldots,\lambda_{m_1}}{\mu_1,\ldots,\mu_{m_2}},\quad
\Lambda'=\binom{b_1',\ldots,b'_{m'_2}}{a'_1,\ldots,a'_{m'_1}}.
\]
Suppose that $\epsilon=+$.
Then we have
\[
\underline\theta_{\bfG'}(\Lambda)
=\theta_0(\Lambda)
=\binom{\tau+m_2,b_1,\ldots,b_{m_2}}{a_1,\ldots,a_{m_1}}
\]
by (\ref{0425}).
We define
\[
\Lambda_0
=\Lambda'^\rmt\smallsetminus\binom{-}{b'_1}
=\binom{a'_1,\ldots,a'_{m'_1}}{b'_2,\ldots,b'_{m'_2}}.
\]
It is clear that $\Lambda'=\underline\theta_{\bfG'}(\Lambda_0)$ from the definition in Subsection~\ref{0412}.
Now we know that $(\Lambda,\Lambda')\in\calb_{\bfG,\bfG'}$, and by Lemma~\ref{0203},
we have the following two situations:
\begin{enumerate}
\item Suppose that $m_1'=m_1$ and $m_2'=m_2+1$, i.e., $|\Lambda^*|=|(\Lambda_0)^*|$ and $|\Lambda_*|=|(\Lambda_0)_*|$.
Now the assumption that $\Lambda'\in\Theta_{\bfG'}(\Lambda)$ implies that
$a_i\geq a'_i>a_{i+1}$ and $b'_i>b_i\geq b'_{i+1}$ for each $i$ by Lemma~\ref{0203}.
The assumption that $\Lambda'\neq\underline\theta_{\bfG'}(\Lambda)$ implies that there exists $1\leq i_0\leq m_1$
such that $a_{i_0}>a'_{i_0}$ or there exists $1\leq j_0\leq m_2$ such that
$b_{j_0}>b'_{j_0+1}$ and hence $n''={\rm rk}(\Lambda_0)<{\rm rk}(\Lambda)=n$ by (\ref{0201}).

\item Suppose that $m_1'=m_1-1$ and $m_2'=m_2$.
Now the assumption that $\Lambda'\in\Theta_{\bfG'}(\Lambda)$ implies that
$a_i>a'_i\geq a_{i+1}$ and $b'_i\geq b_i>b'_{i+1}$ for each $i$ by Lemma~\ref{0203},
and hence $n''={\rm rk}(\Lambda_0)<{\rm rk}(\Lambda)=n$ by Lemma~\ref{0210}.
\end{enumerate}

Next suppose that $\epsilon=-$.
Then we define
\[
\Lambda_0
=\Lambda'^\rmt\smallsetminus\binom{a'_1}{-}
=\binom{a'_2,\ldots,a'_{m'_1}}{b'_1,\ldots,b'_{m'_2}}.
\]
Then $\Lambda'=\underline\theta_{\bfG'}(\Lambda_0)$ and we have the following two situations:
\begin{enumerate}
\item Suppose that $m_1'=m_1+1$ and $m_2'=m_2$.
Now the assumption that $\Lambda'\in\Theta_{\bfG'}(\Lambda)$ implies that
$a'_i>a_i\geq a'_{i+1}$ and $b_i\geq b'_i>b_{i+1}$ for each $i$ by Lemma~\ref{0203}.
The assumption that $\Lambda'\neq\underline\theta_{\bfG'}(\Lambda)$ implies that there exists $i_0$ such that
$a_{i_0}>a'_{i_0+1}$ or $b_{i_0}>b'_{i_0}$ and hence $n''={\rm rk}(\Lambda_0)<{\rm rk}(\Lambda)=n$.

\item Suppose that $m_1'=m_1$ and $m_2'=m_2-1$.
Now the assumption that $\Lambda'\in\Theta_{\bfG'}(\Lambda)$ implies that
$a'_i\geq a_i>a'_{i+1}$ and $b_i>b'_i\geq b'_{i+1}$ for each $i$ by Lemma~\ref{0203},
and hence $n''={\rm rk}(\Lambda_0)<{\rm rk}(\Lambda)=n$ by Lemma~\ref{0210}, again.
\end{enumerate}
\end{proof}

\begin{prop}\label{0305}
Consider the dual pair $(\bfG,\bfG')=(\rmO^\epsilon_{2n},\Sp_{2n'})$ with $2n\leq n'$.
Then $\eta$, $\underline\theta$ and $\overline\theta$ coincide on unipotent characters, i.e.,
\[
\eta(\rho_\Lambda)
=\underline\theta_{\bfG'}(\rho_\Lambda)
=\overline\theta_{\bfG'}(\rho_\Lambda)
\]
for any $\Lambda\in\cals_{\rmO_{2n}^\epsilon}$.
\end{prop}
\begin{proof}
Let $\Lambda\in\cals_{\rmO_{2n}^\epsilon}$.
Now the dual pair is in stable range, so by Proposition~\ref{0504}, we only need to prove the first equality.
We prove the Lemma by induction on $n$.

Suppose that $\epsilon=+$.
First suppose that $n=0$.
Then $\Lambda=\binom{-}{-}$ and $\rho_\Lambda$ is the trivial character $\rmO^+_0(q)$.
It is clear that $\eta(\rho_\Lambda)$ is the trivial character of $\Sp_{2n'}(q)$.
By definition we see that $\underline\theta_{\bfG'}(\Lambda)=\binom{n'}{-}$,
and hence $\underline\theta_{\bfG'}(\rho_\Lambda)=\rho_{\binom{n'}{-}}$ is also the trivial character of $\Sp_{2n'}(q)$.
Therefore, the lemma holds for $n=0$.
Suppose that the lemma is true for any $n''<n$.
Now from Proposition~\ref{0301}, we know that $\eta(\rho_\Lambda)$
is the unique element in $\Theta_{\bfG'}(\rho_\Lambda)$ with rank equal to $n$.
On the other hand, by Lemma~\ref{0302},
any irreducible character $\rho_{\Lambda'}\in\Theta_{\bfG'}(\rho_\Lambda)\smallsetminus\{\underline\theta_{\bfG'}(\rho_\Lambda)\}$
is the image of $\underline\theta_{\bfG'}(\rho_{\Lambda_0})$ for some $\Lambda_0\in\cals_{\rmO_{2n''}^+}$ such that
$n''<n$.
Then by induction hypothesis, $\rho_{\Lambda'}=\eta(\rho_{\Lambda_0})$ and hence
${\rm rk}(\rho_{\Lambda'})=n''<n$ and hence $\rho_{\Lambda'}\neq\eta(\rho_\Lambda)$.
Therefore $\underline\theta_{\bfG'}(\rho_\Lambda)=\eta(\rho_\Lambda)$.

Next suppose that $\epsilon=-$.
Suppose that $n=2$.
Then $\Lambda=\binom{-}{1,0}$ or $\binom{1,0}{-}$.
It is not difficult to see that
$\Theta_{\bfG'}(\binom{-}{1,0})=\{\binom{1,0}{n'}\}$ and
$\Theta_{\bfG'}(\binom{1,0}{-})=\{\binom{-}{n',1,0}\}$ by Lemma~\ref{0203}.
Because now both $\underline\theta_{\bfG'}(\rho_\Lambda)$ and $\eta(\rho_\Lambda)$ are in
$\Theta_{\bfG'}(\rho_\Lambda)$, we have $\underline\theta_{\bfG'}(\rho_\Lambda)=\eta(\rho_\Lambda)$
for $\Lambda=\binom{-}{1,0}$ or $\binom{1,0}{-}$.
The remaining proof is the same as in the case for $\epsilon=+$.
\end{proof}

%% file: sect06.tex
% !TEX root = finite-unipotent.tex

\section{Eta Correspondence and Lusztig Correspondence}
In this section, we consider a reductive dual pair $(\bfG,\bfG')$ of one orthogonal group and one symplectic group.

\subsection{Theta correspondence and Lusztig correspondence}\label{0607}
Let $(\bfG,\bfG')$ be a dual pair of one orthogonal group and one symplectic group.
For a semisimple element $s\in(G^*)^0$,
let groups $G^{(1)}$, $G^{(2)}$, $G^{(3)}$ be defined as in \cite{pan-Lusztig-correspondence}.
Now we have following two cases:
\begin{enumerate}
\item Suppose that the orthogonal group in the dual pair is even.
We have a modified Lusztig correspondence
\[
\Xi_s\colon\cale(G)_s\rightarrow\cale(G^{(1)}\times G^{(2)}\times G^{(3)})_1.
\]
and write $\Xi_s(\rho)=\rho^{(1)}\otimes\rho^{(2)}\otimes\rho^{(3)}$ where
$\rho^{(j)}\in\cale(G^{(j)})_1$ for $j=1,2,3$.

\item Suppose that the orthogonal group in the dual pair is odd.
We have a modified Lusztig correspondence
\[
\Xi'_s\colon\cale(G)_s\rightarrow
\begin{cases}
\cale(G^{(1)}\times G^{(2)}\times G^{(3)})_1, & \text{if $\bfG$ is symplectic};\\
\cale(G^{(1)}\times G^{(2)}\times G^{(3)})_1\times\{\pm1\}, & \text{if $\bfG$ is odd-orthogonal}.
\end{cases}
\]
Write $\Xi'_s(\rho)=\Xi_s(\rho)=\rho^{(1)}\otimes\rho^{(2)}\otimes\rho^{(3)}$ if $\bfG$ is a symplectic group;
and $\Xi'_s(\rho)=\rho^{(1)}\otimes\rho^{(2)}\otimes\rho^{(3)}\otimes\epsilon=\Xi_s(\rho)\otimes\epsilon$
if $\bfG$ is an odd-orthogonal group.
\end{enumerate}

The following proposition is from \cite{pan-Lusztig-correspondence} theorem~6.9 and theorem~7.9:
\begin{prop}\label{0401}
Let $(\bfG,\bfG')$ be a dual pair of one orthogonal group and one symplectic group.
Suppose that $\rho\in\cale(G)_s$ and $\rho'\in\cale(G')_{s'}$.
Keep the above setting.
Then $(\rho,\rho')\in\Theta_{\bfG,\bfG'}$ if and only if the following holds:
\begin{itemize}
\item $\bfG^{(1)}\simeq\bfG'^{(1)}$ and $\rho^{(1)}=\rho'^{(1)}$;

\item $\bfG^{(2)}\simeq\bfG'^{(2)}$ and $\rho^{(2)}=\rho'^{(2)}$;

\item $(\rho^{(3)},\rho'^{(3)})\in\Theta_{\bfG^{(3)},\bfG'^{(3)}}$.
\end{itemize}
i.e., we have the following commutative diagram:
\[
\begin{CD}
\rho @> \Theta_{\bfG'} >> \rho' \\
@V \Xi_s VV @VV \Xi_{s'} V \\
\rho^{(1)}\otimes\rho^{(2)}\otimes\rho^{(3)} @> {\rm id}\otimes{\rm id}\otimes\Theta_{\bfG'^{(3)}} >> \rho'^{(1)}\otimes\rho'^{(2)}\otimes\rho'^{(3)}.
\end{CD}
\]
\end{prop}

Then we can use the Lusztig correspondence to extend the domain of $\underline\theta$ and $\overline\theta$
outside unipotent characters.
First we define
\begin{equation}\label{0402}
\underline\theta_{\bfG'}(\rho)=
\begin{cases}
\Xi_{s'}'^{-1}(\rho^{(1)}\otimes\rho^{(2)}\otimes\underline\theta_{\bfG'^{(3)}}(\rho^{(3)})\otimes\epsilon), & \text{if $\bfG$ is odd-orthogonal};\\
\Xi_{s'}^{-1}(\rho^{(1)}\otimes\rho^{(2)}\otimes\underline\theta_{\bfG'^{(3)}}(\rho^{(3)})), & \text{otherwise}.
\end{cases}
\end{equation}
Note that $\rho^{(3)}$ is unipotent and $\underline\theta_{\bfG'^{(3)}}(\rho^{(3)})$ is defined in (\ref{0427})
and the definition in (\ref{0402}) is to make the $\underline\theta$-correspondence and Lusztig correspondence
commutative, i.e.,
\[
\begin{CD}
\rho @> \underline\theta_{\bfG'} >> \rho' \\
@V \Xi_s VV @VV \Xi_{s'} V \\
\rho^{(1)}\otimes\rho^{(2)}\otimes\rho^{(3)} @> {\rm id}\otimes{\rm id}\otimes\underline\theta_{\bfG'^{(3)}} >> \rho'^{(1)}\otimes\rho'^{(2)}\otimes\rho'^{(3)}.
\end{CD}
\]
Note that when $\bfG$ is an odd-orthogonal group,
the factor $\epsilon$ in (\ref{0402}) is determined by the Witt series of the orthogonal group $\bfG'^{(3)}$.

We define $\overline\theta_{\bfG'}(\rho)$ similarly, i.e., we have a commutative diagram
\[
\begin{CD}
\rho @> \overline\theta_{\bfG'} >> \rho' \\
@V \Xi_s VV @VV \Xi_{s'} V \\
\rho^{(1)}\otimes\rho^{(2)}\otimes\rho^{(3)} @> {\rm id}\otimes{\rm id}\otimes\overline\theta_{\bfG'^{(3)}} >> \rho'^{(1)}\otimes\rho'^{(2)}\otimes\rho'^{(3)}.
\end{CD}
\]
So now we have two mappings $\underline\theta_{\bfG'},\overline\theta_{\bfG'}\colon\cale(G)\rightarrow\cale(G')$,
and then we define
\begin{align*}
\underline\theta_{\bfG,\bfG'} &=\{\,(\rho,\rho')\in\cale(G)\times\cale(G')\mid\rho'=\underline\theta_{\bfG'}(\rho)\,\}; \\
\overline\theta_{\bfG,\bfG'} &=\{\,(\rho,\rho')\in\cale(G)\times\cale(G')\mid\rho'=\overline\theta_{\bfG'}(\rho)\,\}.
\end{align*}

\begin{exam}
Consider the dual pair $(\bfG,\bfG')=(\rmO_1,\Sp_{2n'})$ for some $n'\geq 1$.
\begin{enumerate}
\item Suppose that $\rho=\bf1_{\rmO_1}$.
Then $\Xi'_s(\rho)={\bf1}_{\rmU_0}\otimes{\bf1}_{\Sp_0}\otimes{\bf 1}_{\Sp_0}\otimes1$.
Then we can reduce the correspondence $(\rho,\rho')\in\underline\theta_{\bfG,\bfG'}$ for
$(\bfG,\bfG')=(\rmO_1,\Sp_{2n'})$ via Lusztig correspondence
to the $\underline\theta$-correspondence $(\bf1_{\Sp_0},\rho'^{(3)})$ for $(\bfG^{(3)},\bfG'^{(3)})=(\Sp_0,\rmO^+_{2n'})$.
Because now the dual pair $(\Sp_0,\rmO^+_{2n'})$ is in stable range,
we have $\underline\theta_{\rmO^+_{2n'}}=\overline\theta_{\rmO^+_{2n'}}$.
Note that ${\bf1}_{\Sp_0}=\rho_{\binom{0}{-}}$ and hence
\[
\rho'^{(3)}
=\underline\theta_{\rmO^+_{2n'}}(\rho_{\binom{0}{-}})
=\overline\theta_{\rmO^+_{2n'}}(\rho_{\binom{0}{-}})
=\rho_{\binom{n'}{0}}
={\bf 1}_{\rmO^+_{2n'}}.
\]
Therefore $\rho'$ is the irreducible character in $\cale(\Sp_{2n'}(q))_{s'}$ such that
$\Xi_{s'}(\rho)={\bf1}_{\rmU_0}\otimes{\bf1}_{\Sp_0}\otimes{\bf1}_{\rmO^+_{2n'}}$,
in particular, we know that the degree of $\rho'$ is $\frac{q^{n'}+1}{2}$.

\item Suppose that $\rho=\sgn_{\rmO_1}$.
Then $\Xi'_s(\rho)={\bf1}_{\rmU_0}\otimes{\bf1}_{\Sp_0}\otimes{\bf 1}_{\Sp_0}\otimes(-1)$.
Then we can reduce the $\underline\theta$-correspondence $(\rho,\rho')$ for
$(\bfG,\bfG')=(\rmO_1,\Sp_{2n'})$ via Lusztig correspondence
to the $\underline\theta$-correspondence $(\bf1_{\Sp_0},\rho'^{(3)})$ for $(\bfG^{(3)},\bfG'^{(3)})=(\Sp_0,\rmO^-_{2n'})$.
Now
\[
\rho'^{(3)}
=\underline\theta_{\rmO^-_{2n'}}(\rho_{\binom{0}{-}})
=\overline\theta_{\rmO^-_{2n'}}(\rho_{\binom{0}{-}})
=\rho_{\binom{-}{n',0}}
={\bf 1}_{\rmO^-_{2n'}}.
\]
Therefore $\rho'$ is the irreducible character in $\cale(\Sp_{2n'}(q))_{s'}$ such that
$\Xi_{s'}(\rho)={\bf1}_{\rmU_0}\otimes{\bf1}_{\Sp_0}\otimes{\bf1}_{\rmO^-_{2n'}}$,
in particular, we know that the degree of $\rho'$ is $\frac{q^{n'}-1}{2}$.
\end{enumerate}
\end{exam}

\subsection{Properties of $\underline\theta$- and $\overline\theta$-correspondence}
So now we have the following properties for the correspondences $\underline\theta$ and $\overline\theta$:

\begin{lem}\label{0405}
Let $(\bfG,\bfG')$ be a reductive dual pair of one symplectic group and one orthogonal group,
and let $\rho\in\cale(G)$.
\begin{enumerate}
\item[(i)] If $\underline\theta_{\bfG'}(\rho)$ is defined,
then $\underline\theta_{\bfG'}(\rho)\in\Theta_{\bfG'}(\rho)$.

\item[(ii)] If $\overline\theta_{\bfG'}(\rho)$ is defined,
then $\overline\theta_{\bfG'}(\rho)\in\Theta_{\bfG'}(\rho)$.
\end{enumerate}
\end{lem}
\begin{proof}
For the special case that the orthogonal group in the pair $(G,G')$ is even and $\rho$ is unipotent,
the lemma is obvious from the definitions in Subsection~\ref{0428}.
From the definition of $\underline\theta$ (resp.~$\overline\theta$) and Proposition~\ref{0401},
we know that both $\underline\theta$ (resp.~$\overline\theta$) and $\Theta$ commute with the Lusztig correspondence.
Then we see that the general case follows from the special case immediately.
\end{proof}

\begin{lem}\label{0406}
Both the correspondences $\underline\theta$ and $\overline\theta$ are symmetric, i.e.,
$\underline\theta_{\bfG'}(\rho)=\rho'$ (resp.~$\overline\theta_{\bfG'}(\rho)=\rho'$) if and only if
$\underline\theta_{\bfG}(\rho')=\rho$ (resp.~$\overline\theta_{\bfG}(\rho')=\rho$).
\end{lem}
\begin{proof}
The proof is exactly similar to that of the previous lemma.
\end{proof}

\begin{lem}\label{0408}
Both the correspondences $\underline\theta$ and $\overline\theta$ are one-to-one, i.e.,
for a dual pair $(\bfG,\bfG')$ and $\rho\in\cale(G)$,
there exists at most one $\rho'\in\cale(G')$ such that $\underline\theta_{\bfG'}(\rho)=\rho'$;
similarly, there exists at most one $\rho''\in\cale(G')$ such that $\overline\theta_{\bfG'}(\rho)=\rho''$.
\end{lem}
\begin{proof}
From Lemma~\ref{0430} and the fact that both $\underline\theta$ and $\overline\theta$ are symmetric,
we see that $\underline\theta$ and $\overline\theta$ are one-to-one on unipotent characters when
$(\bfG,\bfG')$ is a dual pair of one symplectic group and one even orthogonal group.

Now for general situation we suppose that $\rho\in\cale(G)_s$ for some semisimple element
$s$ in the connected component of the dual group $G^*$ of $G$.
From \cite{pan-Lusztig-correspondence}, we know that every irreducible character in $\Theta_{\bfG'}(\rho)$
is in the same Lusztig series $\cale(G')_{s'}$ for some unique conjugacy class $(s')$ determined by $s$
(and the Witt series of $\bfG'$).
Then the result for the general situation follows from the fact that both $\underline\theta$ and $\overline\theta$
are compatible with the Lusztig correspondence.
\end{proof}

\begin{prop}\label{0606}
Let $(\bfG,\bfG')$ be a dual pair in stable range.
Then $\underline\theta_{\bfG'}$ and $\overline\theta_{\bfG'}$ coincide.
\end{prop}
\begin{proof}
Let $\rho\in\cale(G)_s$ for some $s$, $\underline\theta_{\bfG'}(\rho)=\rho'$ and $\overline\theta_{\bfG'}(\rho)=\rho''$.
Note that both $\rho',\rho''$ are in $\Theta_{\bfG'}(\rho)$, so they are in the same Lusztig series $\cale(G')_{s'}$ for some $s'$.
Write $\Xi_s(\rho)=\rho^{(1)}\otimes\rho^{(2)}\otimes\rho^{(3)}$,
$\Xi_{s'}(\rho')=\rho'^{(1)}\otimes\rho'^{(2)}\otimes\rho'^{(3)}$,
and $\Xi_{s'}(\rho'')=\rho''^{(1)}\otimes\rho''^{(2)}\otimes\rho''^{(3)}$.
Then we have the following commutative diagram from the definition in Subsection~\ref{0607}:
\[
\begin{CD}
\rho'' @< \overline\theta_{\bfG'} << \rho @> \underline\theta_{\bfG'} >> \rho' \\
@V \Xi_{s'} VV @V \Xi_s VV @VV \Xi_{s'} V \\
\rho''^{(1)}\otimes\rho''^{(2)}\otimes\rho''^{(3)} @< {\rm id}\otimes{\rm id}\otimes\overline\theta_{\bfG'^{(3)}} << \rho^{(1)}\otimes\rho^{(2)}\otimes\rho^{(3)}
@> {\rm id}\otimes{\rm id}\otimes\underline\theta_{\bfG'^{(3)}} >> \rho'^{(1)}\otimes\rho'^{(2)}\otimes\rho'^{(3)}.
\end{CD}
\]
Since the dual pair $(\bfG,\bfG')$ is in stable range,
it is clear that $(\bfG^{(3)},\bfG'^{(3)})$ is also in stable range,
so $\underline\theta_{\bfG'^{(3)}}(\rho^{(3)})=\overline\theta_{\bfG'^{(3)}}(\rho^{(3)})$ by Proposition~\ref{0504}.
Hence $\rho'^{(j)}=\rho''^{(j)}$ for $j=1,2,3$, which implies that $\rho'=\rho''$,
i.e., $\underline\theta_{\bfG'}(\rho)=\overline\theta_{\bfG'}(\rho)$.
\end{proof}

\subsection{First occurrences for $\underline\theta$}
\begin{lem}
Every irreducible character $\rho'\in\cale(G')$ occurs in the
the correspondence $\underline\theta$ (resp.~$\overline\theta$) in a fixed Witt series, i.e.,
there exist a group $\bfG$ in a fixed Witt series and an irreducible character $\rho\in\cale(G)$ such that
$\underline\theta_{\bfG'}(\rho)=\rho'$ (resp.~$\overline\theta_{\bfG'}(\rho)=\rho'$).
\end{lem}
\begin{proof}
If $\rho'$ is unipotent and the dual pair consists of a symplectic group and an even orthogonal group,
then the result is known in Subsection~\ref{0502}.
Because now both $\underline\theta$ and $\overline\theta$ are compatible with the Lusztig correspondence by definition,
the result for general situation follows.
\end{proof}

From the above lemma, now we can define $\underline n_0(\rho')$ and $\overline n_0(\rho')$ for
any $\rho'\in\cale(G')$ as in Subsection~\ref{0502}.
Clearly by Lemma~\ref{0405} for any $\rho'\in\cale(G')$,
we have $n_0(\rho')\leq\underline n_0(\rho')$ and $n_0(\rho')\leq\overline n_0(\rho')$.

\begin{lem}\label{0410}
Let $\rho\in\cale(G)$ and $G'$ varies in a fixed Witt series.
Then
\[
n_0'(\rho)=\underline n_0'(\rho).
\]
\end{lem}
\begin{proof}
If $G'$ is in a Witt series of even orthogonal group and $\rho$ is unipotent,
the result is just Lemma~\ref{0314}.
Now the general case follows from Lemma~\ref{0314} and Proposition~\ref{0401}
immediately.
\end{proof}

By the above lemma we see that the \emph{(non-)preservation principle} for $\Theta$ in \cite{pan-Lusztig-correspondence}
also holds for $\underline\theta$-correspondence:
\begin{enumerate}
\item[(I)] If $G$ is an orthogonal group and $\rho\in\cale(G)$,
then we have
\[
\underline n_0'(\rho)+\underline n_0'(\rho\cdot\sgn)
=2n-\delta(\rho).
\]

\item[(II)] Suppose that $G$ is a symplectic group and $\rho\in\cale(G)$.
\begin{enumerate}
\item If $\bfG'^\epsilon$ for $\epsilon=\pm$ are in two Witt series of even orthogonal groups,
then we have
\[
\underline n_0'^+(\rho)+\underline n_0'^-(\rho)
=2n-\delta(\rho).
\]

\item Suppose that $\bfG'$ is in a Witt series of odd orthogonal groups.
Suppose that $\rho\in\cale(G)_s$ and $\Xi_s(\rho)=\rho^(1)\otimes\rho^{(2)}\otimes\rho^{(3)}$.
Now we know that $\bfG^{(3)}$ is an even orthogonal group.
Let $\rho^\natural\in\cale(G)_s$ such that $\Xi_s(\rho^\natural)=\rho^(1)\otimes\rho^{(2)}\otimes(\rho^{(3)}\cdot\sgn)$.
Then we have
\[
\underline n_0'(\rho)+\underline n_0'(\rho^\natural)
=2n-\delta(\rho).
\]
\end{enumerate}
\end{enumerate}

\subsection{Eta correspondence and Lusztig correspondence}

\begin{lem}\label{0409}
Let $(\bfG,\bfG')=(\bfG_n,\bfG'_{n'})$ be in stable range and $\rho\in\cale(G)$.
Suppose that $\rho'\in\Theta_{\bfG'}(\rho)\smallsetminus\{\underline\theta_{\bfG'}(\rho)\}$.
Then $n_0(\rho')<n$.
\end{lem}
\begin{proof}
Suppose that $\rho\in\cale(G)_s$ and $\rho'\in\cale(G')_{s'}$
for some $s,s'$ in the connected components of dual groups $G^*,G'^*$ of $G,G'$ respectively.
Write
\[
\Xi_s(\rho)=\rho^{(1)}\otimes\rho^{(2)}\otimes\rho^{(3)}\qquad\text{and}\qquad
\Xi_{s'}(\rho')=\rho'^{(1)}\otimes\rho'^{(2)}\otimes\rho'^{(3)}.
\]
The assumption that $\rho'\in\Theta_{\bfG'}(\rho)$ means that $\rho'^{(3)}\in\Theta_{\bfG'^{(3)}}(\rho^{(3)})$
by Proposition~\ref{0401}.
The assumption that $\rho'\neq\underline\theta_{\bfG'}(\rho)$ means that $\rho'^{(3)}\neq\underline\theta_{\bfG'^{(3)}}(\rho^{(3)})$
by the definition of $\underline\theta$.
Then by Lemma~\ref{0302}, there exists $\rho''^{(3)}\in\cale(G''^{(3)})$ such that $\underline\theta_{\bfG'^{(3)}}(\rho''^{(3)})=\rho'^{(3)}$
and $G''^{(3)}$ is of smaller split rank than that of $G^{(3)}$.
Then by Proposition~\ref{0401} again, there exists $\bfG''=\bfG_{n''}$ with $n''<n$ and $\rho''\in\cale(G'')_{s''}$ for some
$s''$ such that $\Xi_{s''}(\rho'')=\rho^{(1)}\otimes\rho^{(2)}\otimes\rho''^{(3)}$,
and hence $\underline\theta_{\bfG'}(\rho'')=\rho'$.
Therefore, $n_0(\rho')<n$.
\end{proof}

\begin{lem}\label{0403}
Consider the dual pair $(\rmO^\epsilon_k,\Sp_{2n'})$ with $k\leq n'$.
Let $\rho\in\cale(\rmO^\epsilon_k(q))_s$ and $\rho'\in\cale(\Sp_{2n'}(q))_{s'}$ for some semisimple elements $s,s'$.
Then $\rho'=\eta(\rho)$ if and only if
\begin{itemize}
\item $G^{(1)}=G'^{(1)}$ and $\rho^{(1)}=\rho'^{(1)}$;

\item $G^{(2)}=G'^{(2)}$ and $\rho^{(2)}=\rho'^{(2)}$;

\item $\rho'^{(3)}=\underline\theta(\rho^{(3)})$.
\end{itemize}
\end{lem}
\begin{proof}
Write $\Xi_s(\rho)=\rho^{(1)}\otimes\rho^{(2)}\otimes\rho^{(3)}$ and
$\Xi_{s'}(\rho')=\rho'^{(1)}\otimes\rho'^{(2)}\otimes\rho'^{(3)}$.
First suppose that $\rho'=\eta(\rho)$.
Then $\rho'\in\Theta_{\bfG'}(\rho)$ and $\rho'$ is of rank $k$ by Proposition~\ref{0301}.
By Proposition~\ref{0401}, we have
$G^{(1)}=G'^{(1)}$ and $\rho^{(1)}=\rho'^{(1)}$;
$G^{(2)}=G'^{(2)}$ and $\rho^{(2)}=\rho'^{(2)}$;
and $\rho'^{(3)}\in\Theta_{\bfG'^{(3)}}(\rho^{(3)})$.
Now if $\rho'^{(3)}\neq\underline\theta_{\bfG'^{(3)}}(\rho^{(3)})$,
then by Lemma~\ref{0302}, there exists $\rho''^{(3)}\in\cale(G''^{(3)})$ such that $\underline\theta_{\bfG'^{(3)}}(\rho''^{(3)})=\rho'^{(3)}$
and $G''^{(3)}$ is in the same Witt series and of smaller split rank than that of $G^{(3)}$.
Then by Proposition~\ref{0401} again, there exists $\bfG''=\rmO_{k'}^\epsilon$ with $k'<k$ and $\rho''\in\cale(G'')_{s''}$ for some
$s''$ such that $\Xi_{s''}(\rho'')=\rho^{(1)}\otimes\rho^{(2)}\otimes\rho''^{(3)}$.
Then, by Proposition~\ref{0401}, we have $\rho'\in\Theta_{\bfG'}(\rho'')$ and hence $\rho'$ is of rank $k'<k$
and we get a contradiction.

Conversely, suppose that
$G^{(1)}=G'^{(1)}$ and $\rho^{(1)}=\rho'^{(1)}$;
$G^{(2)}=G'^{(2)}$ and $\rho^{(2)}=\rho'^{(2)}$;
$\rho'^{(3)}=\underline\theta_{\bfG'^{(3)}}(\rho^{(3)})$.
Then we have $\rho'=\underline\theta_{\bfG'}(\rho)$ by (\ref{0402}).
Now if $\eta(\rho)\neq\rho'=\underline\theta_{\bfG'}(\rho)$,
then by Lemma~\ref{0409} we know that the character $\eta(\rho)$ occurs in the $\Theta$-correspondence for
the dual pair $(\rmO_{k'}^\epsilon(q),\Sp_{2n'}(q))$ for some $k'<k$.
This conflicts with Proposition~\ref{0301}, so we must have $\eta(\rho)=\rho'$.
\end{proof}

\begin{rem}
Consider the dual pair $(\rmO_k^\epsilon,\Sp_{2n'})$ with $k\leq n'$ and $k$ even.
In this case, $\bfG^{(3)}$ is an even orthogonal group and the dual pair $(\bfG^{(3)},\bfG'^{(3)})$ is in stable range.
Then by Lemma~\ref{0305}, we know that $\eta(\rho^{(3)})=\underline\theta_{\bfG'^{(3)}}(\rho^{(3)})$.
Hence the lemma means that $\eta$-correspondence and the Lusztig correspondence are compatible
for a dual pair of an even orthogonal group and a symplectic group and in stable range.
\end{rem}

The following theorem says that the $\underline\theta$-correspondence
(defined for any dual pair of an orthogonal group and a symplectic group)
is an extension of the $\eta$-correspondence (defined only for dual pairs in stable range):

\begin{thm}\label{0404}
Consider the dual pair $(\rmO^\epsilon_k,\Sp_{2n'})$ such that $k\leq n'$.
Then the correspondences $\eta,\underline\theta_{\bfG'},\overline\theta_{\bfG'}$ coincide, i.e.,
\[
\eta(\rho)
=\underline\theta_{\bfG'}(\rho)
=\overline\theta_{\bfG'}(\rho)
\]
for any $\rho\in\cale(\rmO^\epsilon_k(q))$.
\end{thm}
\begin{proof}
Now the dual pair is in stable range, the second equality follows from Proposition~\ref{0606}.
So now we only need to consider the first equality.
If $k$ is even and $\rho$ is unipotent,
the result is Proposition~\ref{0305}.

Now we consider the general situation.
Suppose that $\rho\in\cale(G)_s$ for some $s$, $\underline\theta_{\bfG'}(\rho)=\rho'$ and $\eta(\rho)=\rho''$.
Note that both $\rho',\rho''$ are in $\Theta_{\bfG'}(\rho)$, so they are in the same Lusztig series $\cale(G')_{s'}$ for some $s'$.
Write $\Xi_s(\rho)=\rho^{(1)}\otimes\rho^{(2)}\otimes\rho^{(3)}$,
$\Xi_{s'}(\rho')=\rho'^{(1)}\otimes\rho'^{(2)}\otimes\rho'^{(3)}$,
and $\Xi_{s'}(\rho'')=\rho''^{(1)}\otimes\rho''^{(2)}\otimes\rho''^{(3)}$.
Then we have the following commutative diagram from Lemma~\ref{0403}:
\[
\begin{CD}
\rho'' @< \eta << \rho @> \underline\theta_{\bfG'} >> \rho' \\
@V \Xi_{s'} VV @V \Xi_s VV @VV \Xi_{s'} V \\
\rho''^{(1)}\otimes\rho''^{(2)}\otimes\rho''^{(3)} @< {\rm id}\otimes{\rm id}\otimes\underline\theta_{\bfG'^{(3)}} << \rho^{(1)}\otimes\rho^{(2)}\otimes\rho^{(3)}
@> {\rm id}\otimes{\rm id}\otimes\underline\theta_{\bfG'^{(3)}} >> \rho'^{(1)}\otimes\rho'^{(2)}\otimes\rho'^{(3)}.
\end{CD}
\]
Then we have $\rho'^{(j)}=\rho''^{(j)}$ for $j=1,2,3$.
Therefore $\rho'=\rho''$, i.e., $\underline\theta_{\bfG'}(\rho)=\eta(\rho)$.
\end{proof}

%% file: sect07.tex
% !TEX root = finite-unipotent.tex

\section{Maximal One-to-One Theta Relation}

\subsection{One-to-one theta relation}\label{0407}
We can regard $\Theta$ as a mapping given by
\[
\Theta\colon (\bfG,\bfG')\mapsto \Theta_{\bfG,\bfG'}\subset\cale(G)\times\cale(G').
\]
If $\vartheta$ is another mapping $(\bfG,\bfG')\mapsto\vartheta_{\bfG,\bfG'}\subset\cale(G)\times\cale(G')$ such that
$\vartheta_{\bfG,\bfG'}$ is a subset of $\Theta_{\bfG,\bfG'}$ for each dual pair $(\bfG,\bfG')$,
then $\vartheta$ is called a \emph{sub-relation} of $\Theta$.
For a sub-relation $\vartheta$ and $\rho\in\cale(G)$,
we define
\[
\vartheta_{\bfG'}(\rho)=\{\,\rho'\in\cale(G')\mid(\rho,\rho')\in\vartheta_{\bfG,\bfG'}\,\}.
\]
Among the set of all sub-relations of $\Theta$,
we can give a partial ordering by inclusion.
Moreover precisely, for two sub-relations $\vartheta^1,\vartheta^2$ of $\Theta$ we say that
$\vartheta^1\subseteq\vartheta^2$ if for each dual pair $(\bfG,\bfG')$,
we have $\vartheta^1_{\bfG,\bfG'}\subseteq\vartheta^2_{\bfG,\bfG'}$ as subsets of $\cale(G)\times\cale(G')$.

Now we have several definitions:
\begin{itemize}
\item A sub-relation $\vartheta$ of $\Theta$ is called \emph{semi-persistent} (on unipotent characters)
if it satisfies the following conditions:
\begin{enumerate}
\item if either 
\begin{enumerate} 
\item $(\bfG,\bfG')=(\Sp_{2n},\rmO^+_{2n'})$ and $\Lambda\in\cals_{\bfG}$ with ${\rm def}(\Lambda)=4d+1$; or

\item $(\bfG,\bfG')=(\rmO^+_{2n'},\Sp_{2n})$ and $\Lambda\in\cals_{\bfG}$ with ${\rm def}(\Lambda)=4d$,
\end{enumerate}
then $\vartheta_{\bfG'}(\rho_\Lambda)\neq\emptyset$ for any $n'\geq n-2d$;

\item if either 
\begin{enumerate}
\item $(\bfG,\bfG')=(\Sp_{2n},\rmO^-_{2n'})$ and $\Lambda\in\cals_{\bfG}$ with ${\rm def}(\Lambda)=4d+1$; or

\item $(\bfG,\bfG')=(\rmO^-_{2n'},\Sp_{2n})$ and $\Lambda\in\cals_{\bfG}$ with ${\rm def}(\Lambda)=4d+2$,
\end{enumerate}
then $\vartheta_{\bfG'}(\rho_\Lambda)\neq\emptyset$ for any $n'\geq n+2d+1$.
\end{enumerate}

\item A sub-relation $\vartheta$ of $\Theta$ is called \emph{symmetric} if
for each dual pair $(\bfG,\bfG')$, any $\rho\in\cale(G)$ and $\rho'\in\cale(G')$,
we have $\rho'\in\vartheta_{\bfG'}(\rho)$ if and only if $\rho\in\vartheta_{\bfG}(\rho')$.

\item A sub-relation $\vartheta$ of $\Theta$ is said to be \emph{compatible with the Lusztig correspondence},
if for each dual pair $(\bfG,\bfG')$, any $\rho\in\cale(G)_s$ and $\rho'\in\cale(G')_{s'}$
for some semisimple elements $s,s'$,
there exist choices of Lusztig correspondences $\Xi_s,\Xi_{s'}$ such that the following diagram
\[
\begin{CD}
\rho @> \vartheta_{\bfG'} >> \rho' \\
@V \Xi_s VV @VV \Xi_{s'} V \\
\rho^{(1)}\otimes\rho^{(2)}\otimes\rho^{(3)} @> {\rm id}\otimes{\rm id}\otimes\vartheta_{\bfG'^{(3)}} >> \rho'^{(1)}\otimes\rho'^{(2)}\otimes\rho'^{(3)}.
\end{CD}
\]
is commutative.
\end{itemize}

A sub-relation $\vartheta$ of $\Theta$ is called a \emph{theta-relation} if it is
semi-persistent, symmetric and compatible with the Lusztig correspondence.

A sub-relation $\vartheta$ of $\Theta$ is called \emph{one-to-one} if for each dual pair $(\bfG,\bfG')$
and each $\rho\in\cale(G)$,
there exists at most one $\rho'\in\cale(G')$ such that $(\rho,\rho')\in\vartheta_{\bfG,\bfG'}$.
If $\vartheta$ is a one-to-one sub-relation of $\Theta$, then for a dual pair $(\bfG,\bfG')$ and $\rho\in\cale(G)$
we know that either $\vartheta_{\bfG'}(\rho)=\emptyset$ or $\vartheta_{\bfG'}(\rho)=\{\rho'\}$ for some
$\rho'\in\cale(G')$.
For the latter case, we may just abuse the notation a little by writing $\vartheta_{\bfG'}(\rho)=\rho'$.

\begin{prop}\label{0604}
Both $\underline\theta$ and $\overline\theta$ are one-to-one theta-relations.
\end{prop}
\begin{proof}
By definition, it is clear that both $\underline\theta$ and $\overline\theta$ are sub-relations of $\Theta$.
When $(\bfG,\bfG')$ consists of a symplectic group and an even orthogonal group,
both $\underline\theta$ and $\overline\theta$ are clearly semi-persistent from the definitions in
Subsection~\ref{0428}.
We know that both $\underline\theta$ and $\overline\theta$ are symmetric by Lemma~\ref{0406},
and they are compatible with the Lusztig correspondence by the definition in (\ref{0402}).
Therefore, both $\underline\theta$ and $\overline\theta$ are theta-relations.
Finally, both $\underline\theta$ and $\overline\theta$ are one-to-one by Lemma~\ref{0408}.
\end{proof}

\subsection{Maximal one-to-one theta-relation}
The $\eta$-correspondence is only defined for dual pairs in stable range.
Now we know that both $\underline\theta$ and $\overline\theta$ are
one-to-one theta-relations which extend $\eta$ to general dual pairs.
The following proposition means that both $\underline\theta$ and $\overline\theta$
can not be extended any more if we require the extension to be a one-to-one theta-relation. 

\begin{prop}\label{0605}
No one-to-one theta-relation can be contained properly in another one-to-one theta-relation.
\end{prop}
\begin{proof}
Let $\vartheta$ be a one-to-one theta-relation.
Suppose that $\vartheta'$ is a theta-relation which properly contains $\vartheta$.
So we need to show that $\vartheta'$ is not one-to-one.
Let $(\bfG,\bfG')$ be a dual pair, $\rho\in\cale(G)$ and $\rho'\in\cale(G')$
such that $(\rho,\rho')\in\vartheta'_{\bfG,\bfG'}$ and $(\rho,\rho')\not\in\vartheta_{\bfG,\bfG'}$.
If $\vartheta_{\bfG'}(\rho)$ is defined, i.e., $(\rho,\rho'')\in\vartheta_{\bfG,\bfG'}\subset\vartheta_{\bfG,\bfG'}'$ for some $\rho''\in\cale(G')$,
then $\vartheta'$ is not one-to-one.
So we may assume that $\vartheta_{\bfG'}(\rho)$ is not defined.
Since both $\vartheta'$ and $\vartheta$ are compatible with Lusztig correspondence,
we have a commutative diagram
\[
\begin{CD}
\rho @> \vartheta'_{\bfG'} >> \rho' \\
@V \Xi_s VV @VV \Xi_{s'} V \\
\rho^{(1)}\otimes\rho^{(2)}\otimes\rho^{(3)} @> {\rm id}\otimes{\rm id}\otimes\vartheta'_{\bfG'^{(3)}} >> \rho'^{(1)}\otimes\rho'^{(2)}\otimes\rho'^{(3)}.
\end{CD}
\]
This means that $\vartheta'_{\bfG'^{(3)}}(\rho^{(3)})=\rho'^{(3)}$ and $\vartheta_{\bfG'^{(3)}}(\rho^{(3)})$ is undefined.
Now $(\bfG^{(3)},\bfG'^{(3)})$ is a dual pair of one symplectic group and one even orthogonal group,
and both $\rho^{(3)},\rho'^{(3)}$ are unipotent.
So we write $\rho^{(3)}=\rho_\Lambda$ and $\rho'^{(3)}=\rho_{\Lambda'}$ for some
$\Lambda\in\cals_{\bfG^{(3)}}$ and $\Lambda'\in\cals_{\bfG'^{(3)}}$.

Suppose that $(\bfG^{(3)},\bfG'^{(3)})=(\Sp_{2n},\rmO^+_{2n'})$ for some $n,n'$,
$\Lambda\in\cals_{\bfG^{(3)}}$ with ${\rm def}(\Lambda)=4d+1$ for some $d\in\bbZ$.
Then ${\rm def}(\Lambda')=-4d=4(-d)$ by Proposition~\ref{0204}.
Because now $\vartheta$ is semi-persistent and $\vartheta_{\bfG'^{(3)}}(\rho_\Lambda)$ is not defined,
we must have $n'<n-2d$.
So we have $n>n'-(-2d)$, and this means that $\vartheta_{\bfG^{(3)}}(\rho_{\Lambda'})$ is defined.
Let $\Lambda''=\vartheta_{\bfG^{(3)}}(\Lambda')\in\cals_{\bfG^{(3)}}$.
Now $(\rho_{\Lambda'},\rho_{\Lambda''})\in\vartheta_{\bfG'^{(3)},\bfG^{(3)}}\subseteq\vartheta'_{\bfG'^{(3)},\bfG^{(3)}}$ implies
that $(\rho_{\Lambda''},\rho_{\Lambda'})\in\vartheta'_{\bfG^{(3)},\bfG'^{(3)}}$ since $\vartheta'$ is symmetric.
Moreover, $(\rho_\Lambda,\rho_{\Lambda'})$ is in $\vartheta'_{\bfG^{(3)},\bfG'^{(3)}}$ by our assumption.
However, $\Lambda$ and $\Lambda''$ are not equal
because $\vartheta_{\bfG'^{(3)}}(\Lambda)$ is not defined and $\vartheta_{\bfG'^{(3)}}(\Lambda'')=\Lambda'$
by the symmetricity of $\vartheta$.
So we conclude that $\vartheta'$ is not one-to-one.

The proof for other cases (i.e., $(\bfG^{(3)},\bfG'^{(3)})=(\Sp_{2n},\rmO^-_{2n'})$, $(\rmO^+_{2n},\Sp_{2n'})$, or $(\rmO^-_{2n},\Sp_{2n'})$)
are similar.
\end{proof}

\begin{cor}\label{0701}
Both $\underline\theta$ and $\overline\theta$ are maximal one-to-one theta-relation, i.e.,
they are not properly contained in any other one-to-one theta-relation.
\end{cor}
\begin{proof}
This follows from Proposition~\ref{0604} and Proposition~\ref{0605} immediately.
\end{proof}